\documentclass[a4paper]{amsart} 
\usepackage{amssymb,times, amscd,amsmath,amsthm, xypic}
\usepackage{tikz}
\usepackage{tikz-cd}
\usepackage[all]{xy}
\usepackage{geometry}
\usepackage{mathrsfs}

\begin{document} 

\numberwithin{equation}{section}
\newtheorem{thm}[equation]{Theorem}
\newtheorem{pro}[equation]{Proposition}
\newtheorem{prob}[equation]{Problem}
\newtheorem{qu}[equation]{Question}
\newtheorem{cor}[equation]{Corollary}
\newtheorem{con}[equation]{Conjecture}
\newtheorem{lem}[equation]{Lemma}
\theoremstyle{definition}
\newtheorem{ex}[equation]{Example}
\newtheorem{defn}[equation]{Definition}
\newtheorem{observation}[equation]{Observation}
\newtheorem{rem}[equation]{Remark}

\renewcommand{\rmdefault}{ptm}
\newcommand{\pdarrow}[2]{\ar@<0.5ex>[d]^-{#1} \ar@<-0.5ex>[d]_-{#2}}

\def \calB{\mathcal B}
\def\frak{\mathfrak}
\def\alp{\alpha}
\def\be{\beta}
\def\jeden{1\hskip-3.5pt1}
\def\om{\omega}
\def\bigstar{\mathbf{\star}}
\def\ep{\epsilon}
\def\vep{\varepsilon}
\def\Om{\Omega}
\def\la{\lambda}
\def\La{\Lambda}
\def\si{\sigma}
\def\Si{\Sigma}
\def\Cal{\mathcal}
\def\m {\mathcal}
\def\mC {\mathcal C}
\def\ga{\gamma}
\def\Ga{\Gamma}
\def\de{\delta}
\def\De{\Delta}
\def\bF{\mathbb{F}}
\def\bH{\mathbb H}
\def\bPH{\mathbb {PH}}
\def \bB{\mathbb B}
\def \bA{\mathbb A}
\def \bC{\mathbb C}
\def \bZ{\mathbb Z}
\def \bOB{\mathbb {OB}}
\def \bM{\mathbb M}
\def \bOM{\mathbb {OM}}
\def \mA{\mathcal A}
\def \mB{\mathscr B}
\def \mC{\mathscr C}
\def \mD{\mathscr D}
\def \mS{\mathscr S}
\def \mR{\mathcal R}
\def \mH{\mathcal H}
\def \mM{\mathcal M}
\def \mV{\mathcal V}
\def \mTOP{\mathcal {TOP}}
\def \mAB{\mathcal {AB}}
\def \bI{\mathbb I}
\def \bK{\mathbb K}
\def \bG{\mathbf G}
\def \bL{\mathbb L}
\def\bN{\mathbb N}
\def\bR{\mathbb R}
\def\bP{\mathbb P}
\def\bZ{\mathbb Z}
\def\bC{\mathbb  C}
\def \bQ{\mathbb Q}
\def\op{\operatorname}

\newcommand{\maru}[1]{\ooalign{
\hfil\resizebox{.8\width}{\height}{#1}\hfil
\crcr
\raise.1ex\hbox{\large$\bigcirc$}}}

\newcommand{\Maru}[1]{\ooalign{
\hfil\resizebox{.6\width}{\height}{#1}\hfil
\crcr
\raise.1ex\hbox{\LARGE$\bigcirc$}}}

\newcommand{\MMaru}[1]{\ooalign{
\hfil\resizebox{.5\width}{\height}{#1}\hfil
\crcr
\raise.1ex\hbox{\Huge$\bigcirc$}}}

\newcommand{\MMMaru}[1]{\ooalign{
\hfil\resizebox{.4\width}{\height}{#1}\hfil
\crcr
\raise.1ex\hbox{\Huge$\bigcirc$}}}

\newcommand{\MMMMaru}[1]{\ooalign{
\hfil\resizebox{.3\width}{\height}{#1}\hfil
\crcr
\raise.1ex\hbox{\Huge$\bigcirc$}}}

\newcommand{\MMMMMaru}[1]{\ooalign{
\hfil\resizebox{.2\width}{\height}{#1}\hfil
\crcr
\raise.1ex\hbox{\Huge$\bigcirc$}}}

\title[A co-operational bivariant theory
derived from cohomology operations]
{A co-operational bivariant theory\\
derived from cohomology operations
}


\thanks {
\noindent
\emph{keywords} : bivariant theory, operational bivariant theory, cohomology operation \\
\emph{Mathematics Subject Classification 2000}: 55N35, 55S99, 14F99}

\author{Shoji Yokura}

\date{}
\address{Graduate School of Science and Engineering, Kagoshima University, 1-21-35 Korimoto, Kagoshima, 890-0065, Japan}
\email{yokura@sci.kagoshima-u.ac.jp}

\begin{abstract} A co-operational bivariant theory of a contravariant functor is a ``dual" version of Fulton--MacPherson's operational bivariant theory of a covariant functor. In this paper we show that these co-operational and operational bivariant theories can be both described as the end of the corresponding bivariant functors. For a given contravariant functor we define a ``restricted" generalized cohomology operation for continuous maps having sections, using cohomology operations. This ``restricted" generalized cohomology operation is motivated by Quillen's Steenrod power operation. Using this generalized cohomology operation, we define \emph{a co-operational bivariant theory derived from cohomology operations}, 
  which is a subtheory of the co-operational bivariant theory associated to the contravariant functor.
\end{abstract}

\maketitle

\section{Introduction}

In \cite[Part I]{FM} (also see \cite[Chapter 17 Bivariant Intersection Theory]{Ful}) William Fulton and Robert MacPherson have introduced Bivariant Theory unifying covariant functors and contravariant functors. A bivariant theory $\mathbb B(X \xrightarrow f Y)$ is a graded abelian group (sometimes a non-graded abelian group or just a set) assigned to a morphism $f:X \to Y$, \emph{not to an object $X$},  in a category $\mathcal C$. Then for a morphism $a_X:X \to pt$ to a point pt, $\mathbb B_*(X):= \mathbb B^*(X \xrightarrow {a_X} pt)$ becomes a covariant functor and for the identity morphism $\op{id}_X:X \to X$, $\mathbb B^*(X):= \mathbb B^*(X \xrightarrow {\op{id}_X} X)$ becomes a contravariant functor. That is why it is called a \emph{\underline{bi}}variant theory. A transformation of two bivariant theories $\ga: \mathbb B \to \mathbb B'$ is called a Grothendieck transformation, which specializes to natural transformations $\ga_*:\mathbb B_* \to \mathbb B'_*$ and $\ga^*:\mathbb B^* \to \mathbb B'^*$ of the associated covariant and contravariant functors, respectively.

As a typical example of a bivariant theory, Fulton and MacPherson have constructed a bivariant homology theory $\mathbb Bh(X \to Y)$ from the usual cohomology theory $h^*$ such that the associated contravariant functor $\mathbb Bh^*$ is isomorphic to the original contravariant functor $h^*$ (see \cite[\S 3 Topological Theories]{FM}). A key fact in the construction of $\mathbb Bh(X \to Y)$ is that the cohomology theory $h^*$ is \emph{multiplicative} \cite[\S 3.1]{FM}. Applying this construction to a natural transformation $t: h^* \to h'^*$ of multiplicative cohomology theories $h^*, h'^*$, they obtain a Grothendieck transformation
$\ga: \mathbb Bh \to \mathbb Bh'$
such that the associated natural transformation $\ga^*:\mathbb Bh^* \to \mathbb Bh'^*$ is the same as the original natural transformation $t:h^* \to h'^*$. For example, they applied it to the Chern character $ch: K \to H^*$.

A characteristic class $c\ell$ of complex vector bundles is a \emph{natural transformation} 
$c\ell: \op{Vect} \to H^*$
from the \emph{set\footnote{In fact, $\op{Vect}(X)$ is a commutative \emph{rig} (``ring" with ``n" deleted) (e.g., see \cite{Schau}, \cite{Law1, Law2}), i.e., a commutative monoid with respect to $\oplus$ and a commutative monoid with respect to $\otimes$, thus a commutative ring without the requirement that every element has an additive inverse.}-valued} contravariant functor $\op{Vect}$ of the isomorphism classes of complex vector bundles to the usual cohomology functor $H^*$ \emph{considered just as a set-valued functor} 
and 
it is well-known that any characteristic class $c\ell$ is a formal power series $p(c_1, c_2, \cdots)$ on Chern classes $c_i$ where $c_i(E) \in H^{2i}(X)$ for a complex vector bundle over a space $X$ (e.g. see \cite[p.321]{Mac2}). Thus we have
$$c\ell=p(c_1, c_2, \cdots): \op{Vect} \to H^*.$$

Since the above functor $\op{Vect}$ is 
\emph{not a multiplicative cohomology theory}, we \emph{cannot} apply the above Fulton--MacPherson's construction $\ga: \mathbb Bh \to \mathbb Bh'$ 
to $c\ell:\op{Vect} \to H^*$. 
In order to be able to obtain a Grothendieck transformation from a natural transformation between two contravariant functors, in \cite{Yo-MRR} we have introduced a \emph{co-operational bivariant theory}, which can be considered as a \emph{``dual" version} of Fulton--MacPherson's operational bivariant theory \cite[\S 8.1]{FM} which is defined for a given covariant functor.

Given a contravariant functor $F^*$, our co-operational bivariant theory, denoted $\mathbb B^{coop}F^*(X \xrightarrow f Y)$, is constructed in a manner analogous to the definition (see \cite[\S 8.1]{FM}) of Fulton--MacPherson's operational bivariant theory $\mathbb B^{op}F_*(X \xrightarrow f Y)$\footnote{In \cite[\S 8.1]{FM} it is simply denoted by $C^*(X \xrightarrow f Y)$.} defined for a covariant functor $F_*$. 

The definition of $\mathbb B^{coop}F^*(X \xrightarrow f Y)$ is the following (for more details, see \S 3). For a map $f:X \to Y$, an element $c \in \mathbb B^{coop}F^*(X \xrightarrow f Y)$ 
is defined to be
\begin{equation}\label{c_g}
c= \Bigl (c_g:F^*(X') \to F^*(Y') \Bigr )_{g \in \mathcal C/Y}
\end{equation}
Here $\mathcal C/Y$ is the slice category\footnote{$Obj(\mathcal C/Y)$ consists of all morphisms to a given object $Y$. For two objects $g':X' \to Y$ and $g:X \to Y$ in $\mathcal C/Y$, a morphism $h:g' \to g$ is a morphism $h:X' \to X$ such that the following diagram commutes $\xymatrix@R=15pt@C=10pt{X' \ar[dr]_{g'} \ar[rr]^h && X \ar[dl]^g \\ & Y}$.} over $Y$ and the map $c_g:F^*(X') \to F^*(Y')$ is such that for the following fiber square on the left, the following diagram on the right commutes:
\begin{equation}\label{f.s.}
\xymatrix@R=15pt@C=20pt{
X'' \ar[d]_{f''} \ar[r]^{h'} & X' \ar[d]_{f'} \ar[r]^{g'} & X \ar[d]^f\\
Y'' \ar[r]_h & Y' \ar[r]_g & Y,
} \qquad \quad 
\xymatrix@R=15pt@C=20pt
{F^*(X'') \ar[d]_{c_{g \circ h}} && F^*(X')\ar[ll]_{h'^*} \ar[d]^{c_g}\\
F^*(Y'') && F^*(Y') \ar[ll]^{h^*}.
}
\end{equation}

It turns out that our co-operational bivariant theory and Fulton--MacPherson's operational bivariant theory can be both expressed as \emph{the end of the corresponding bivariant functors} exactly in the same way (see Theorem \ref{thm-coop-op} ).

A more refined co-operational bivariant theory, which is a ``dual" version of Fulton--MacPherson's refined operational bivariant theory \cite[\S 8.3]{FM}, will be treated in \cite{Yo2}.

The well-known cohomology operation $\theta$ for a cohomology theory $h^*$ is by the definition (e.g., see \cite{Steenrod1, Steenrod2, Steenrod3, Mo-Tan, nLab2, nLab}) a collection of maps $\theta_X:h^*(X) \to h^*(X)$, one for each space $X$, such that for every continuous map $f:X \to Y$ the following diagram commutes:
\begin{equation}\label{co-op}
\xymatrix@R=20pt@C=20pt
{ h^*(X) \ar[d]_{\theta_X} & h^*(Y) \ar[l]_{f^*} \ar[d]^{\theta_Y}\\
h^*(X) & h^*(Y). \ar[l]^{f^*}
}
\end{equation}
Category-Functor theoretically, it can be called \emph{a natural (self-)transformation from $h^*$ to itself}. 
So, any natural (self-)transformation $\theta: F^* \to F^*$ from a contravariant functor $F^*$ to itself shall be still called \emph{a cohomology operation of $F^*$}, instead of ``a contravariant functor operation". 

Let $coop F^*$ be the set of all cohomology operations of $F^*$. It turns out (Propostion \ref{coopF^*}) that we have 
$$coopF^*= \bB^{coop}F^*(pt \xrightarrow {\op{id}_{pt}} pt)$$ where $pt$ is the terminal object. 
For the map $p_X: X \to pt$, using the following bivariant pullback
$$p_X^*: \bB^{coop}F^*(pt \xrightarrow {\op{id}_{pt}} pt) \to   \bB^{coop}F^*(X \xrightarrow {\op{id}_X} X).$$
we define
$$coopF^*(X \xrightarrow {\op{id}_X} X):= \op{Image}(p_X^*) \subset \bB^{coop}F^*(X \xrightarrow {\op{id}_X} X).$$
For a cohomology operation $\theta \in coopF^*$, the image $p_X^*(\theta)$ is explicitly described by
$$p_X^*(\theta) = \Bigl ( (p_X^*(\theta))_g=\theta_{X'}:F^*(X') \to F^*(X') \Bigr )_{g \in \mathcal C/X}$$
where we consider the following trivial fiber square
$$\xymatrix@R=15pt@C=20pt{
X' \ar[d]_{\op{id}_{X'}} \ar[r]^g & X \ar[d]^{\op{id}_X}\\
X' \ar[r]_g & X
}
$$
Here we emphasize that since $ (p_X^*(\theta))_g=\theta_{X'}$ for $g:X' \to X (\in \mathcal C/X)$, we do have \emph{$(p_X^*(\theta))_g = (p_X^*(\theta))_{g'}$ for any $g, g' \in hom(X',X)$ }.

Suggested by the definition of $\mathbb B^{coop}F^*(X \xrightarrow f Y)$ and the definition of cohomology operation of $F^*$, it is natural or reasonable to define an operation $\widehat \theta$ assigning a map $\widehat \theta_f: F^*(X) \to F^*(Y)$ for \emph{every} continuous map $f:X \to Y$ such that for the following fiber square on the left, the following diagram on the right commutes:
\begin{equation}\label{co-op2}
\xymatrix@R=15pt@C=20pt{
X' \ar[d]_{f'} \ar[r]^{g'} & X \ar[d]^f\\
Y' \ar[r]_g & Y,
} \qquad \quad \xymatrix@R=15pt@C=20pt
{F^*(X') \ar[d]_{\widehat \theta_{f'}} && F^*(X)\ar[ll]_{g'^*} \ar[d]^{\widehat \theta_f}\\
F^*(Y') && F^*(Y) \ar[ll]^{g^*}.
}
\end{equation}
Such an operation $\widehat \theta$ shall be called \emph{a generalized cohomology operation} of $F^*$ since it becomes a cohomology operation when it is restricted to \emph{every} identity map $\op{id}_X:X \to X$. Indeed, for the following fiber square on the left, the following diagram on the right commutes:
\begin{equation*}
\xymatrix@R=20pt@C=20pt{
X' \ar[d]_{\op{id}_{X'}} \ar[r]^g & X \ar[d]^{\op{id}_X}\\
X' \ar[r]_g & X
}\qquad \quad \xymatrix
{F^*(X') \ar[d]_{\theta_{X'} :=\widehat \theta_{\op{id}_{X'}}} && F^*(X)\ar[ll]_{g^*} \ar[d]^{\widehat \theta_{\op{id}_X}=:\theta_X}\\
F^*(X') && F^*(X) \ar[ll]^{g^*}.
}
\end{equation*}
As far as we know, it seems that no work has been done for such a generalized cohomology operation even for the well-known cohomology theory $H^*$ and the complex cobordism $MU^*$, although a lot of work have been done for the cohomology operations of them.

A cohomology operation $\theta =(\theta_X)_{X \in \mathcal C}$ is indexed by all the objects. On the other hand, a co-operational bivariant element $c = (c_g)_{g \in \mathcal C/Y} \in \bB^{coop}F^*(X \xrightarrow f Y)$ \emph{indexed by the objects of the slice categiry $\mathcal C/Y$} is defined for \emph{the maps $f':X' \to Y'$ pulled back from a fixed map $f:X \to Y$ by all maps $g:Y' \to Y$} in the slice category $\mathcal C/Y$ over the target space $Y$ of the map $f:X \to Y$. In this sense, such a cohomology operation is called a \emph{``restricted" generalized cohomology operation}.

For the complex cobordism $U^*$, Daniel Quillen (\cite[Proposition 3.17, p.43]{Qui}, cf. \cite[\S 4 Steenrod Operations]{EG}) considers 
$$P: U^{-2q}(X) \to U^{-2qk}(B \times X)
$$
which is defined by the composition\footnote{This kind of Steenrod power operation can be defined for \emph{any} map $f:Z \to X$, since $P:U^{-2q}(X) \to U^{-2qk}(Z)$ can be defined as the composition
$U^{-2q}(X) \xrightarrow {f^*} U^{-2q}(Z) \xrightarrow {\theta_Z} U^{-2qk}(Z)$.}
$U^{-2q}(X) \xrightarrow {\pi^*} U^{-2q}(B \times X) \xrightarrow {\theta_{B \times X}} U^{-2qk}(B \times X)$
where $\pi: B \times X \to X$ is the projection to the second factor and $\theta_{B \times X}$ is the Steenrod operation.
He calls the above $P$ the ``Steenrod $k$-th power operation".
This assignment $P$ is \emph{not} a ``restricted" generalized cohomology operation for such a map $\pi:B \times X \to X$.  

In this paper, motivated by the description of Quillen's Steenrod $k$-th power operation, i.e., the composition of the pullback and a cohomology operation, we will construct \emph{a co-operational bivariant theory derived from  cohomology operations $\delta\bB^{coop}F^*(X \xrightarrow f Y)$}, using the above $coopF^*(X \xrightarrow {\op{id}_X} X)$ and the pullback $s^*$ for sections $s:Y \to X$, for a map $f:X \to Y$ having sections $s:Y \to X$.


The paper is organized as follows. In \S 2 we give a quick recall of Fulton--MacPherson's bivariant theory \cite{FM}. In \S 3 we give a quick recall of a co-operational bivariant theory \cite{Yo-MRR}. In \S 4 we give an ``unified" description of our co-operational bivariant theory and Fulton--MacPherson's operational bivariant theory, using the notion of \emph{\underline{end} of a bivariant functor}. 
In \S 5 we discuss ``restricted" generalized cohomology operations for special maps called ``\emph{sectional maps}", which are, as mentioned above, motivated by Steenrod $k$-th power operation. 
In \S 6 we construct a co-operational bivariant theory derived from the data of cohomology operations. 

\section{Fulton--MacPherson's bivariant theories}\label{BT}
We give a very quick recall of Fulton--MacPherson's bivariant theory \cite {FM}.

Let $\mathscr C$ be a category which has a final object $pt$ and on which the fiber product or fiber square is well-defined. $\mathscr C$ is equipped with the following classes:
\begin{enumerate}
\item a class $\mathcal C$ of maps, called ``confined maps" (e.g., proper maps, in algebraic geometry), 
\item a class $\mathcal Ind$ of commutative diagrams, called ``independent squares" (e.g., fiber square, ``Tor-independent" square, in algebraic geometry)
\begin{equation}\label{ind-sq}
\xymatrix@R=15pt@C=15pt{
X' \ar[d]_{f'} \ar[r]^{g'} & X \ar[d]^f\\
Y' \ar[r]_g & Y.
}
\end{equation}
\end{enumerate}
These are required to satisfy the following conditions:
\begin{enumerate}
\item The class $\mathcal C$ of confined maps are \emph{closed under composition and base change, and contain all the identity maps},

\item if the two inside squares in  
$$
\xymatrix{ X'' \ar[d]_{f''}  \ar[r]^{h'} & X' \ar[d]_{f'} \ar[r]^{g'} & X \ar[d]^ f\\
Y'' \ar[r]_{h} & Y' \ar[r]_{g} & Y, 
} 
\quad \quad \qquad \text{or} \qquad \quad \quad 
\xymatrix@R=15pt@C=20pt
{X' \ar[d]_{f'} \ar[r]^{h''} & X \ar[d]^f\\
Y' \ar[d]_{g'} \ar[r]^{h'} & Y \ar[d]^g\\
Z' \ar[r]_{h} & Z 
} 
$$
are independent\footnote{We note that give an independent square 
its transpose is \emph{not necessarily} independent. For example, let us consider the category of topological spaces and continuous maps. Let \emph{any} map be confined, and we let a fiber square (\ref{fib.squ})
to be \emph{independent only if $g$ is proper} (hence $g'$ is also proper). Then its transpose is \emph{not independent unless 
$f$ is proper}. (Note that the pullback of a proper map by any continuous map is proper, because ``proper" is equivalent to ``universally closed", i.e., the pullback by any map is closed.)
}, then the outside square is also independent,

\item any square of the following forms are independent:
$$
\xymatrix@R=15pt@C=15pt{X \ar[d]_{f} \ar[r]^{\op {id}_X}&  X \ar[d]^f & & X \ar[d]_{\op {id}_X} \ar[r]^f & Y \ar[d]^{\op {id}_Y} \\
Y \ar[r]_{\op {id}_X}  & Y && X \ar[r]_f & Y}
$$
where $f:X \to Y$ is \emph{any} morphism. 
\item In any independent square (\ref{ind-sq}), if $f$ (resp., $g$) is confined, then $f'$ (resp., $g'$) is also confined.
\end{enumerate}

\begin{def}\label{BivariantTheory}
A \emph{bivariant theory} $\mathbb B$ on a category $\mathscr C$ with values in the category of graded abelian groups\footnote{As we mentioned in Introduction, instead of graded abelian groups, we consider also sets, e.g., such as the set of complex structures and the set of Spin structures (see \cite[\S 4.3.2]{FM}),  and categories, e.g., such as the derived (triangulated) category of $f$-perfect complexes (see \cite[\S 7.1 Grothendieck duality]{FM}) as well.}
 is an assignment to each morphism
$ X  \xrightarrow{f} Y$
in the category $\mathscr C$ a graded\footnote{The grading is sometimes ignored. In this case we can consider that the grading is only $0$, i.e, $\bB(X \to Y) = \bB^0(X \to Y)$. For an example for such a case, see \cite[\S 6.1.2 Definition of $\mathbb F$]{FM} where $\mathbb F(X \to Y) = \mathbb F^0(X \to Y)$.}
abelian group $$\bB(X  \xrightarrow{f} Y)$$
which is equipped with the following three basic operations. The $i$-th component of $\bB(X  \xrightarrow{f} Y)$, $i \in \mathbb Z$, is denoted by $\bB^i(X  \xrightarrow{f} Y)$.
\begin{enumerate}
\item {\bf Product}: For morphisms $f: X \to Y$ and $g: Y
\to Z$, the product 
$$\bullet: \bB^i( X  \xrightarrow{f}  Y) \otimes \bB^j( Y  \xrightarrow{g}  Z) \to
\bB^{i+j}( X  \xrightarrow{g\circ f}  Z).$$
\item {\bf Pushforward}: For morphisms $f: X \to Y$
and $g: Y \to Z$ with $f$ \emph {confined}, 

the pushforward 
$f_*: \bB^i( X  \xrightarrow{g\circ f} Z) \to \bB^i( Y  \xrightarrow{g}  Z).$
\item {\bf Pullback} : For an \emph{independent} square \qquad 
$\xymatrix@R=15pt@C=15pt{
X' \ar[d]_{f'} \ar[r]^{g'} & X \ar[d]^f\\
Y' \ar[r]_g & Y.
}$

the pullback 
$g^* : \bB^i( X  \xrightarrow{f} Y) \to \bB^i( X'  \xrightarrow{f'} Y'). $
\end{enumerate}

An element $\alp \in \bB(X \xrightarrow f Y)$ is sometimes expressed as follows:
\[
\xymatrix
{
X \ar[rr]_f^{\maru{$\alp$}} && Y
} 
\]

These three 
operations are required to satisfy the following seven compatibility 
axioms (\cite [Part I, \S 2.2]{FM}):

\begin{enumerate}
\item[($A_1$)] {\bf Product is associative}: for 
\xymatrix
{
X \ar[r]_f^{\maru{$\alp$}} & Y \ar[r]_g^{\maru{$\be$}} & Z \ar[r]_h^{\maru{$\ga$}} & Z
}
$$(\alp \bullet\be) \bullet \ga = \alp \bullet (\be \bullet \ga).$$
\item[($A_2$)] {\bf Pushforward is functorial}:for 
\xymatrix
{
X  \ar@/^10pt/[rrr]^{\maru{$\alp$}}  \ar[r]_f & Y \ar[r]_g  & Z \ar[r]_h  & W
}
with confined $f, g$, 
$$(g\circ f)_* \alp = g_*(f_*\alp).$$
\item[($A_3$)] {\bf Pullback is functorial}: given independent squares
$
\xymatrix@R=15pt@C=20pt{
X'' \ar[d]_{f''} \ar[r]^{h'} & X' \ar[d]_{f'} \ar[r]^{g'} & X \ar[d]^f\\
Y'' \ar[r]_h & Y' \ar[r]_g & Y,
} 
$
$$(g \circ h)^* = h^* \circ g^*.$$
\item[($A_{12}$)] {\bf Product and pushforward commute}: $f_*(\alp \bullet\be)  = f_*\alp \bullet \be$ \\

for 
\xymatrix
{
X \ar@/^10pt/[rr]^{\maru{$\alp$}}  \ar[r]_f & Y \ar[r]_g  & Z \ar[r]_h^{\maru{$\be$}} & W
}
with confined $f$, 
\item[($A_{13}$)] {\bf Product and pullback commute}: \quad $h^*(\alp \bullet\be)  = {h'}^*\alp \bullet h^*\be$ \quad for independent squares
 \[
\xymatrix@R=15pt@C=20pt
{X' \ar[rrr]^{h''} \ar[d]_{f'}
 &&& X \ar[d]_f^{\maru{$\alp$}}   \\
Y' \ar[d]_{g'}\ar[rrr]^{h'} &&& Y \ar[d]_g^{\maru{$\be$}} \\
Z' \ar[rrr]_h &&& Z
} 
\]
\item[($A_{23}$)] 
{\bf Pushforward and pullback commute}: \quad $f'_*(h^*\alp)  = h^*(f_*\alp)$ \quad for independent squares with $f$ confined 
\[
\xymatrix@R=15pt@C=20pt{
X' \ar[rrr]^{h''} \ar[d]_{f'}  &&& X \ar[d]_f \ar@/^11pt/[dd]^{\maru{$\alp$}} \\
Y' \ar[d]_{g'} \ar[rrr]^{h'} &&& Y \ar[d]_g\\
Z' \ar[rrr]_h &&& Z
} 
\]
\item[($A_{123}$)] {\bf Projection formula}: \, $g'_*(g^*\alp \bullet \be)  = \alp \bullet g_*\be$ \quad for an independent square with $g$ confined 
\[
\xymatrix{
X' \ar[r]^{g'} \ar[d]^{f'}_{\Maru{$g^*\alp$}} 
& X \ar[d]_f^{\maru{$\alp$}}  \\
Y' \ar@/_18pt/[rrrrr]^{\maru{$\beta$}} \ar[r]^g & Y \ar[rrrr]_h^{\quad \Maru{$g_*\be$}} &&&& Z
} 
\]
\end{enumerate}
We also require the theory $\bB$ to have multiplicative units:
\begin{enumerate}
\item[({\bf Units})] For all $X \in \mathscr C$, there is an element $1_X \in \bB^0( X  \xrightarrow{\op {id}_X} X)$ such that $\alp \bullet 1_X = \alp$ for all morphisms $W \to X$ and all $\alp \in \bB(W \to X)$, and such that $1_X \bullet \beta = \beta $ for all morphisms $X \to Y$ and all $\beta \in \bB(X \to Y)$, and such that $g^*1_X = 1_{X'}$ for all $g: X' \to X$.
\end{enumerate}
\end{def}

A bivariant theory unifies both a covariant theory and a contravariant theory in the following sense:
For a bivariant theory $\bB$, its associated covariant functors and contravariant functors are defined as follows:
\begin{enumerate}
\item $\bB_*(X):= \bB(X \xrightarrow {a_X} pt)$ is covariant for confined morphisms and the grading is given by $\bB_i(X):= \bB^{-i}(X \xrightarrow {a_X} pt)$.
\item $\bB^*(X) := \bB(X  \xrightarrow{\op{id}_X}  X)$ is contravariant for all morphisms and the grading is given by $\bB^j(X):= \bB^j(X   \xrightarrow{\op{id}_X}  X)$.
\end{enumerate}

A typical example of a bivariant theory is the bivariant homology theory $\bH(X \xrightarrow f Y)$ 
constructed from the singular cohomology theory $H^*(-)$, which unifies the Borel--Moore homology $H_*^{BM}(X):=\bH^{-*}(X \to pt)$ and the singular cohomology $H^*(X):=\bH^*(X \xrightarrow {\op{id}_X} X)$. Here the underlying category $\mathscr C$ is the category of spaces embeddable as closed subspaces of some Euclidean spaces $\mathbb R^n$ and continuous maps between them (see \cite[\S 3 Topological Theories]{FM}).
More generally, Fulton--MacPherson's (general) bivariant homology theory 
$$h^*(X \to Y)$$
(here, using their notation) 
is constructed from \emph{a multiplicative cohomology theory $h^*(-)$} \cite[\S 3.1]{FM}.
Here the cohomology theory $h^*$ is either ordinary or generalized. 
A cohomology theory $h^*$ is called \emph{multiplicative} if for pairs $(X,A)$ and $(Y,B)$, there is a graded pairing (exterior product)
$$h^i(X,A) \times h^j(Y,B) \xrightarrow {\times} h^{i+j}(X \times Y, X \times B \sqcup A \times Y)$$
such that it is associative and graded commutative, i.e., $\alp \times \be = (-1)^{i+j} \be \times \alp$.
A typical example of a multiplicative ordinary cohomology theory is the singular cohomology theory. The topological complex $K$-theory $K(-)$ and complex cobordism theory $\Omega^*(-)$ are  multiplicative generalized  cohomology theories. 
\begin{defn}\cite[\S 2.7 Grothendieck transformation]{FM}\label{groth}
Let $\bB$ and $\bB'$ be two bivariant theories on a category $\mathscr C$. 
A {\it Grothendieck transformation} from $\bB$ to $\bB'$, $\ga : \bB \to \bB'$
is a collection of homomorphisms
$\bB(X \to Y) \to \bB'(X \to Y)$
for a morphism $X \to Y$ in the category $\mathscr C$, which preserves the above three basic operations: 
\begin{enumerate}
\item $\ga (\alp \bullet_{\bB} \be) = \ga (\alp) \bullet _{\bB'} \ga (\be)$, 
\item $\ga(f_{*}\alp) = f_*\ga (\alp)$, and 
\item $\ga (g^* \alp) = g^* \ga (\alp)$. 
\end{enumerate}
\end{defn}

A Grothendieck transformation $\ga: \bB \to \bB'$ induces natural transformations $\ga_*: \bB_* \to \bB_*'$ and $\ga^*: \bB^* \to {\bB'}^*$.
 
\section{Co-operational bivariant theory}\label{coop}

We make a quick recall of our co-operational bivariant theory \cite{Yo-MRR}. 

\begin{defn}\cite[Theorem 4.1]{Yo-MRR}\label{coop-biv} (Co-operational bivariant theory)
Let $F^*(-)$ be a contravariant functor with values in the category of graded abelian groups, defined on a category $\mathscr C$ with all fiber squares independent.
For each  
map $f:X \to Y$ an element
\begin{equation*}
c \in \mathbb B^{coop}F^i(X \xrightarrow f Y)
\end{equation*}
is defined to be a collection of 
maps\footnote{In our original paper \cite[Theorem 4.1]{Yo-MRR}, we required $c_g$ to be a homomorphism, as in Fulton--MacPherson's operational bivariant theory \cite[Part I, \S 8 Operational Theories]{FM}.}
\begin{equation}\label{coll}
c_g: F^m(X') \to F^{m+i}(Y')
\end{equation}
for all $m \in \mathbb Z$, all $g:Y' \to Y$  and the fiber square \, \, \, 
\begin{equation}
\xymatrix@R=15pt@C=15pt{
X' \ar[d]_{f'} \ar[r]^{g'} & X \ar[d]^f\\
Y' \ar[r]_g & Y
}
\end{equation}
And these maps $\{c_g\}$ are
required to be compatible with pullback, i.e., for the following fiber diagram on the left, the following diagram on the right commutes:
\begin{equation}\label{coop-cd}
\xymatrix@R=15pt@C=15pt
{ X'' \ar[d]_{f''}  \ar[r]^{h'} & X' \ar[d]_{f'} \ar[r]^{g'} & X \ar[d]^ f\\
Y'' \ar[r]_{h} & Y' \ar[r]_{g} & Y, 
} \qquad 
\xymatrix@R=15pt@C=15pt
{F^m(X'') \ar[d]_{c_{g \circ h}} && F^m(X')\ar[ll]_{h'^*} \ar[d]^{c_g}\\
F^{m+i}(Y'') && F^{m+i}(Y') \ar[ll]^{h^*}.
}
\end{equation}
We define the following bivariant-theoretic operations:
\begin{enumerate}
\item {\bf Product}: The product
$$\bullet: \mathbb B^{coop}F^i(X \xrightarrow f Y) \otimes \mathbb B^{coop}F^j(Y \xrightarrow g Z) \to \mathbb B^{coop}F^{i+j}(X \xrightarrow {g \circ f} Z)$$
is defined by, for $c \in \mathbb B^{coop}F^i(X \xrightarrow f Y) $ and $d \in \mathbb B^{coop}F^j(Y \xrightarrow g Z) $,
$$(c \bullet d)_h:= d_h \circ c_{h'}: F^m(X') \xrightarrow {c_{h'}}  F^{m+i}(Y') \xrightarrow {d_h} F^{m+i+j}(Z').$$
Here we consider the following fiber squares:
\begin{equation}\label{cd-1}
\xymatrix@R=15pt@C=20pt
{X' \ar[d]_{f'} \ar[r]^{h''} & X \ar[d]^f\\
Y' \ar[d]_{g'} \ar[r]^{h'} & Y \ar[d]^g\\
Z' \ar[r]_{h} & X 
} 
\end{equation}
\item {\bf Pushforward}: For $X \xrightarrow f Y \xrightarrow g Z$ 
$$f_*: \mathbb B^{coop}F^i(X \xrightarrow {g \circ f} Z) \to \mathbb B^{coop}F^i(Y \xrightarrow g Z)$$
is defined by, for $c \in \mathbb B^{coop}F^i(X \xrightarrow {g \circ f} Z)$, 
$$(f_*c)_h := c_h \circ f'^*: F^m(Y') \xrightarrow{f'^*} F^m(X')  \xrightarrow{c_h} F^{m+i}(Z').$$
Here we use the above commutative diagram (\ref{cd-1}).
\item {\bf Pullback}: For a fiber square 
$$
\xymatrix@R=15pt@C=15pt{
X' \ar[d]_{f'} \ar[r]^{g'} & X \ar[d]^f\\
Y' \ar[r]_g & Y
}
$$
$$g^*: \mathbb B^{coop}F^i(X \xrightarrow f Y ) \to \mathbb B^{coop}F^i(X' \xrightarrow {f'} Y' )$$
is defined by, for $c \in \mathbb B^{coop}F^i(X \xrightarrow f Y) $,
$$(g^*c)_h := c_{g \circ h}: F^m(X'') \to F^{m+i}(Y'').$$
Here we use the commutative diagram on the left of (\ref{coop-cd}).
\end{enumerate} 
\end{defn}
\begin{thm} \cite[Proof of Theorem 4.1]{Yo-MRR}.
$\mathbb B^{coop}F^*$
is a bivariant theory, i.e., satisfies the seven axioms of bivariant theory\footnote{Although we do not require $c_g$ to be a homomorphism in the above Definition \ref{coop-biv}, this statement is still correct.}.
\end{thm}
\begin{defn} \cite[Definition 4.22]{Yo-MRR}\label{4-defn-bcoop}
Let $T:F^* \to G^*$ be a natural transformation between two contravariant functors.
For $f:X \to Y$,  
an element
$$c \in  \mathbb B^{coop}_T F^i(X \xrightarrow f Y) \quad (\subset \mathbb B^{coop}F^i(X \xrightarrow f Y) )$$
is defined to satisfy
that there exists an element  
$$c^T \in  \mathbb B^{coop} G^i(X \to Y)$$
such that
the following diagram commutes:for $m \in \mathbb Z$ and for $g:Y' \to Y$
\begin{equation}\label{cd-7}
\xymatrix@R=10pt@C=15pt
{F^m(X') \ar[d]_{c_g} \ar[rr]^T && G^m(X') \ar[d]^{(c^T)_g}\\
F^{m+i}(Y') \ar[rr]_T && G^{m+i}(Y').
}
\end{equation}
\end{defn}
\begin{thm}{\cite[Theorem 4.26]{Yo-MRR}}.
The above $\mathbb B^{coop}_T F^*$ is a bivariant theory. Namely, the following product, pushforward and pullback are all well-defined.
\begin{enumerate}
\item {\bf Product}:\quad $\bullet: \mathbb B^{coop}_TF^i(X \xrightarrow f Y) \otimes \mathbb B^{coop}_TF^j(Y \xrightarrow g Z) \to \mathbb B^{coop}_TF^{i+j}(X \xrightarrow {g \circ f} Z)$,
\item {\bf Pushforward}: \quad $f_*: \mathbb B^{coop}_TF^*(X \xrightarrow {g \circ f} Z) \to \mathbb B^{coop}_TF^*(Y \xrightarrow g Z)$,
\item {\bf Pullback}: For a fiber square 
$$
\xymatrix@R=15pt@C=15pt{
X' \ar[d]_{f'} \ar[r]^{g'} & X \ar[d]^f\\
Y' \ar[r]_g & Y
}
$$
$$g^*: \mathbb B^{coop}_TF^*(X \xrightarrow f Y ) \to \mathbb B^{coop}_TF^*(X' \xrightarrow {f'} Y' ).$$
\end{enumerate} 
\end{thm}
\begin{cor}{\cite[Corollary 4.35]{Yo-MRR}}.\label{cor-nat} Let $T:F^* \to G^*$ be a natural transformation between two contravariant functors. For $\mathbb B^{coop}_T F^*$ we assume that bivariant elements $c^T$ in Definition \ref{4-defn-bcoop} satisfy the following three properties (if $T$ is surjective, they are automatically satisfied, as pointed out in \cite[Remark 4.34]{Yo-MRR}):
\begin{enumerate}
\item $(c \bullet d)^T = c^T \bullet d^T$,
\item $(f_*c)^T =f_*(c^T)$, 
\item $(g^*c)^T =g^*(c^T)$. 
\end{enumerate}
\noindent 
Then we have a Grothendieck transformation
$$\ga_T: \mathbb B^{coop}_T F^* \to  \mathbb B^{coop} G^*$$
which is a collection of 
$\ga_T: \mathbb B^{coop}_T F^*(X \xrightarrow f Y)  \to  \mathbb B^{coop} G^*(X \xrightarrow f Y)$
defined by, for each $c \in \mathbb B^{coop}_T F^*(X \to Y)$, 
$$\ga_T(c) := c^T.$$
Furthermore $\ga_T(c) = c^T$ means that for $m \in \mathbb Z$ and for $g:Y' \to Y$ we have the commutative diagram (\ref{cd-7}), i.e., \begin{equation*}
\xymatrix@R=10pt@C=15pt
{F^m(X') \ar[d]_{c_g} \ar[rr]^T && G^m(X') \ar[d]^{(c^T)_g}\\
F^{m+i}(Y') \ar[rr]_T && G^{m+i}(Y').
}
\end{equation*}
\end{cor}
\begin{rem}
For the sake of later presentation, here we recall some of Fulton--MacPherson's operational bivariant theory.
We emphasize that there are, in fact, \emph{two} theories, which are 
\begin{enumerate}
\item \emph{a foundational one} defined in \cite[\S 8.1, pp.88--90]{FM} and 
\item \emph{a more refined one} defined in \cite[\S 8.3, p.91]{FM} with more requirements on the foundational one. 
\end{enumerate}
We recall the foundational one, since our co-operational bivariant theory is a ``dual" version of this foundational one.

Let $F_*$ be a homology theory, i.e., a covariant functor on the subcategory of confined maps with values in graded abelian groups, thus we have the pushforward homomorphism $f_*:F_*(X) \to F_*(Y)$ for a confined map $f:X \to Y$.  Then its \emph{foundational operational bivariant theory} is denoted by $\bB^{op}F_*(X \xrightarrow f Y)$\footnote{In \cite[\S 8 Operational Theories]{FM}, it is denoted by $C^i(X \xrightarrow f Y)$.} and an element $c \in \bB^{op}F_*^i(X \xrightarrow f Y)$ is defined as follows:
For a map $f:X \to Y$, an element $c \in \bB^{op}F^i_*(X \xrightarrow f Y)$ is defined to be a collection of
homomorphisms
\begin{equation*}
c_g:F_m(Y') \to F_{m-i}(X')
\end{equation*}
for all $m \in \mathbb Z$, all $g:Y' \to Y$ and the fiber square \, \, \, 
$$
\xymatrix@R=13pt@C=15pt{
X' \ar[d]_{f'} \ar[r]^{g'} & X \ar[d]^f\\
Y' \ar[r]_g & Y
}
$$

\noindent
These homomorphisms $c_g$ are
required to be \emph{compatible with pushforwards} (for confined maps),  i.e., for a fiber diagram on the left (where $k$ is \emph{confined, thus $k'$ is confined as well}) the diagram on the right commutes:
\begin{equation}\label{op-cd}
\xymatrix@R=15pt@C=15pt
{ X'' \ar[d]_{f''}  \ar[r]^{h'} & X' \ar[d]_{f'} \ar[r]^{g'} & X \ar[d]^ f\\
Y'' \ar[r]_{h} & Y' \ar[r]_{g} & Y, 
} \qquad 
\xymatrix@R=15pt@C=15pt
{F_*(X'') \ar[rr]^{h'_*}&& F_*(X') \\
F_*(Y'') \ar[u]^{c_{g \circ h}} \ar[rr]_{h_*}&& F_*(Y') \ar[u]_{c_g}.
}
\end{equation}
\end{rem}
\section{Operational and co-operational bivariant theories can be described as the ``end" of bivariant functors}\label{end}
In this section we show that Fulton--MacPherson's operational bivariant theory and our co-coperational bivariant theory can be described as \emph{the end of a bivariant functor}, which was originally introduced by Nobuo Yoneda \cite{Yoneda} (e.g., also see Mac Lane's book \cite{Maclane}, and \cite{Lor, Lor-2}).

First we need to quickly recall the notion of wedge and end of a bivariant functor.
\begin{defn}[Wedge] Let $\mC$ and $\mD$ be categories and let $B: \mC^{op} \times \mC \to \mD$ be a bivariant functor. 
A wedge of the functor $B$ is a pair $(E, \omega)$ consisting of 
\begin{enumerate}
\item an object $E \in \mD$ and 
\item a collection $\omega =\{\omega_X\}$ of morphisms $\omega_X:E \to B(X,X)$ for all the objects $X \in \mC$ such that 
for every morphism $f:X \to Y$ we have the commutative diagram
\begin{equation}
\xymatrix@R=15pt@C=20pt{
& E \ar[dl]_{\omega_X} \ar[dr]^{\omega_Y}\\
B(X,X) \ar[dr]_{B(\op{id}_X, f)\,\, \, \, }  && B(Y,Y) \ar[dl]^{\,\, B(f, \op{id}_Y)} \\
& B(X,Y).
}
\end{equation}
\end{enumerate}
\end{defn}
\begin{defn}[End]
If a wedge $(E,\omega)$ of a bivariant functor $B: \mC^{op} \times \mC \to \mD$ is called an \emph{end} if it is universal, namely, for any wedge $(E', \omega')$, i.e., for any collection $\omega'$ of morphisms $\{\omega'_X: E' \to B(X,X) \, | \, X \in Obj(\mC)\}$
such that if the outer square of the diagram (\ref{e-bxx}) below is commutative, then there is a unique map $u:E' \to E$ such that the diagram (\ref{e-bxx}) below commutes
\begin{equation}\label{e-bxx}
\xymatrix@R=20pt@C=20pt{
& E' \ar[dl]_{\omega'_X} \ar[dr]^{\omega'_Y} \ar[d]_u\\
B(X,X) \ar[dr]_{B(\op{id}_X, f) \,\, \, \, } & E \ar[l]_{\omega_X} \ar[r]^{\omega_Y} &  B(Y,Y) \ar[dl]^{\,\, B(f, \op{id}_Y)} \\
& B(X,Y)
}
\end{equation}
The end of the functor $B$ is denoted by
$$\int_{X \in \mC} B(X, X).$$
\end{defn}
\begin{thm}(e.g., see \cite[(1.27), p.15]{Lor-2}, cf. his book \cite{Lor} as well) Let the situation be as above. We have
$$\int_{X \in \mC} B(X,X) = \Bigl \{ (c_X)_{X \in \mC} \in \prod_{X \in \mC} B(X,X) \, \, \Bigl | \, \, B(\op{id}_X, f)(c_X) = B(f, \op{id}_Y)(c_Y), \forall f:X \to Y \Bigr \}.$$
Which is the equalizer of the following parallel pair $(B(\op{id}_X, f), B(f, \op{id}_Y))$
\begin{equation}
\xymatrix{
\op{Eq} \Bigl \{ \displaystyle \prod_{X \in \mC} B(X,X) \ar@<0.8ex>[rr]^{B(\op{id}_X, f)} \ar@<-0.8ex>[rr]_{B(f, \op{id}_Y)}  && \displaystyle \prod_{f:X \to Y} B(X,Y)  \Bigr \}
}
\end{equation}
\end{thm}
Now 
we can identify Fulton--MacPherson's operational bivariant theory and our co-coperational bivariant theory as the end of the corresponding bivariant functors as follows, where we use the following notation for the sake of presentation:for a fiber square
\begin{equation}\label{coop-fibsq}
\xymatrix@R=15pt@C=15pt{
X' \ar[d]_{f'} \ar[r]^{g'} & X \ar[d]^f\\
Y' \ar[r]_g & Y
}
\end{equation}
the map $f'$ pullbacked by $g$ of $f$ shall be denotedd by $g^*f$. Then the source space $X$ of a map $f:X \to Y$ is denoted by $dom(f)$ and so on.
\begin{thm}\label{thm-coop-op} Let $\mS$ be the category of sets and maps.
\begin{enumerate}
\item Let $F^*:\mC^{op} \to \mD$ be a contravariant functor to the category $\mD$ of rigs\footnote{In geometry and topology, a contravariant functor such as the cohomology theory and the complex $K$-theory is valued in the category of rings, but in the case of the contravariant functor $\op{Vect}$ is valued in the category of rigs. So we consider this}. For a morphsism $f:X \to Y$ 
we define $BF^*(f):(\mC/Y)^{op} \times (\mC/Y)  \to \mS$ by
\begin{equation*}\label{BF^*}
BF^*(f)(k,g):= hom_{\mS}\Bigl (F^*(dom(g^*f)), F^*(dom(k) \Bigr ),
\end{equation*}
which is a bivariant functor. Then we have 
\begin{equation}\label{BF^*-end}
\int_{k \in \mC/Y} BF^*(f)(k,k) = \int_{k:\mC/Y} hom_{\mS}\Bigl (F^*(dom(k^*f)), F^*(dom(k) \Bigr ) = \bB^{coop}F^*(X \xrightarrow f Y).
\end{equation}
\item Let $F_*:\mC \to \mD$ be a covariant functor to the category $\mD$ of graded abelian groups. For a morphism $f:X \to Y$ we define $BF_*(f):(\mC/Y)^{op} \times (\mC/Y)  \to \mS$ by
\begin{equation*}\label{BF_*}
BF_*(f)(k,g):= hom_{\mS}\Bigl (F_*(dom(k)), F_*(dom(g^*f) \Bigr ),
\end{equation*}
which is a bivariant functor. Then we have 
\begin{equation*}\label{BF_*-end}
\int_{k \in \mC/Y} BF_*(f)(k,k) = \int_{k:\mC/Y} hom_{\mS}\Bigl (F_*(dom(k)), F_*(dom(k^*f) \Bigr ) = \bB^{op}F_*(X \xrightarrow f Y).
\end{equation*}
\end{enumerate}
\end{thm}

\begin{proof} 
We show only (1) in detail and the proof of (2) is left for the reader, since it similar to (1). 

First we show that $BF^*(f):(\mC/Y)^{op} \times (\mC/Y)  \to \mS$ is a bivariant functor. For the sake of simplicity, we just write $BF^*$ deleting $(f)$. For that, first we show the contravariance, i.e., for a fixed $s \in \mC/Y$, $BF^*(-,s):(\mC/Y)^{op} \to \mD$ is a contravariant functor, i.e., for $\xymatrix{\ell \ar[r]^m & k \ar[r]^h & g}$ in $\mC/Y$, we have
\begin{equation}\label{BF^*-contra}
BF^*(h \circ m, \op{id}_s)= BF^*(m, \op{id}_s) \circ BF^*(h, \op{id}_s),
\end{equation}
i.e., for the diagrams on the left, we have the commutative diagram on the right:
$$\xymatrix
{Y''' \ar[dr]_{\ell} \ar[r]^m & Y'' \ar[d]_k \ar[r]^h & Y' \ar[dl]^g\\
& Y,
}\qquad 
\xymatrix@R=30pt@C=1pt
{BF^*(\ell,s) && BF^*(k,s) \ar[ll]_{BF^*(m, \op{id}_s)} \\
& BF^*(g,s) \ar[ul]^{BF^*(h \circ m, \op{id}_s) \quad } \ar[ur]_{\quad BF^*(h, \op{id}_s)}
}
$$
The diagram on the right above is the following (from now on we drop $\mS$, for simplicity):
\begin{equation}\label{BIGCD-1}
\xymatrix@R=30pt@C=-22pt{
hom(F^*(dom(s^*f)), F^*(Y''')) && hom(F^*(dom(s^*f)), F^*(Y''))\ar[ll]_{\quad hom(\op{id}_{F^*(dom(s^*f))}, m^*)} \\
& hom(F^*(dom(s^*f)), F^*(Y'))\ar[ul]^{hom(\op{id}_{F^*(dom(s^*f))}, (h \circ m)^*) \qquad \quad  } \ar[ur]_{\qquad hom({F^*(dom(s^*f))}, h^*)}
}
\end{equation}
Here $h^*:F^*(Y') \to F^*(Y'')$ and $m^*:F^*(Y'') \to F^*(Y''')$ and we have $(h \circ m)^*=m^* \circ h^*$. Thus the above diagram (\ref{BIGCD-1}) is commutative, i.e., (\ref{BF^*-contra}) holds.

Next we show the covariance of $BF^*:(\mC/Y)^{op} \times (\mC/Y)  \to \mS$, i.e., for a fixed $s \in \mC/Y$, $BF^*(s,-):\mC/Y) \to \mD$ is a covariant functor, i.e., for $\xymatrix{\ell \ar[r]^m & k \ar[r]^h & g}$ in $\mC/Y$, we have
\begin{equation}\label{BF^*-co}
BF^*(\op{id}_s, h \circ m)= BF^*(\op{id}_s, h) \circ BF^*(\op{id}_s, m),
\end{equation}
i.e., we have the commutative diagram:
$$
\xymatrix@R=30pt@C=5pt
{BF^*(s, \ell ) \ar[dr]_{BF^*(\op{id}_s, h \circ m) \quad } \ar[rr]^{BF^*(\op{id}_s, m)}&& BF^*(s,k) \ar[dl]^{\quad BF^*(\op{id}_s), h)}\\
& BF^*(s, g), 
}
$$
which is the following:
\begin{equation*}\label{BIGCD-2}
\xymatrix@R=30pt@C=-18pt{
hom( F^*(X'''), F^*(dom(s))) \ar[dr]_{hom((h' \circ m')^*, \op{id}_{F^*(dom(s))}) \quad \qquad } \ar[rr]^{\quad hom(m'^*, \op{id}_{F^*(dom(s))})} && hom( F^*(X''), F^*(dom(s))) \ar[dl]^{\qquad hom(h'^*, \op{id}_{F^*(dom(s))})} \\
& hom( F^*(X'), F^*(dom(s))) }
\end{equation*}
Here we consider the following fiber squares:
$$\xymatrix@R=15pt@C=15pt
{ X''' \ar[d]_{f'''}  \ar[r]^{m'} & X''\ar[d]_{f''} \ar[r]^{h'} & X' \ar[d]^{f'}\\
Y''' \ar[r]_{m} & Y'' \ar[r]_{h} & Y', 
}
$$ 
For $\alp \in hom( F^*(X'''), F^*(dom(s)))$ we have
$hom(m'^*, \op{id}_{F^*(dom(s))})(\alp)= \alp \circ m'^*.$
Hence 
\begin{align*}
hom(h'^*, \op{id}_{F^*(dom(s))}) & \circ hom(m'^*, \op{id}_{F^*(dom(s))}) \\
& = hom(m'^* \circ h'^*, \op{id}_{F^*(dom(s))})\\
& = hom((h' \circ m')^*, \op{id}_{F^*(dom(s))})
\end{align*}
Hence we get (\ref{BF^*-co}). 

Now we show (\ref{BF^*-end}). 
For $f:X \to Y$, $g:Y' \to Y \in \mC/Y$ and $k:Y'' \to Y \in \mC/Y$ with $h:Y'' \to Y'$ (note $k=g \circ h$), we have the following fiber squares:
$$
\xymatrix@R=15pt@C=15pt
{ X'' \ar[d]_{f''}  \ar[r]^{h'} & X' \ar[d]_{f'} \ar[r]^{g'} & X \ar[d]^ f\\
Y'' \ar[r]_{h} & Y' \ar[r]_{g} & Y. 
}
$$  
$dom(k^*f)=dom(f'')=X''$, $dom(g^*f)=dom(f')=X'$, $dom(k)=Y''$ and $dom(g)=Y'$, thus
$$BF^*(k, k) = hom(F^*(dom(k^*f)), F^*(dom(k))) = hom(F^*(X''), F^*(Y'')),$$
$$BF^*(k,g):= hom(F^*(dom(g^*f)), F^*(dom(k)))= hom(F^*(X'), F^*(Y'')),$$
$$BF^*(g,g):= hom(F^*(dom(g^*f)), F^*(dom(g)))= hom(F^*(X'), F^*(Y')).$$
Thus the following cospan on the left is the cospan on the right:
$$
\xymatrix@R=20pt@C=-20pt
{BF^*(k,k) \ar[dr]_{BF^*(\op{id}_k, h)} && BF^*(g,g) \ar[dl]^{BF^*(h, \op{id}_g)} \\
& BF^*(k,g),
}\quad 
\xymatrix@R=20pt@C=-45pt
{hom(F^*(X''), F^*(Y'')) \ar[dr]_{hom(h'^*, \op{id}_{F^*(Y'')})\quad } && hom(F^*(X'), F^*(Y')).\ar[dl]^{\quad hom(\op{id}_{F^*(X')}, h^*)} \\
& hom(F^*(X'), F^*(Y''))
}
$$
Thus the parallel pair 
$$
\xymatrix{
\displaystyle \prod_{k \in \mC/Y} BF^*(k,k) \ar@<1ex>[rrr]^{BF^*(\op{id}_k, h)} \ar@<-1ex>[rrr]_{BF^*(h, \op{id}_g)}  &&& \displaystyle \prod_{h:k \to g} BF^*(k,g)
}
$$ 
is the following
$$
\xymatrix{
\displaystyle \prod_{k \in \mC/Y} hom\Bigl (F^*(X''), F^*(Y'') \Bigr ) \ar@<1ex>[rrr]^{hom(h'^*, \op{id}_{F^*(Y'')})} \ar@<-1ex>[rrr]_{hom(\op{id}_{F^*(X')}, h^*)}  &&& \displaystyle \prod_{h:k \to g} hom\Bigl (F^*(X'), F^*(Y'') \Bigr ) 
}
$$
whose equalizer is exactly the following (note that $k= g \circ h$, thus $c_k=c_{g \circ h}$)  
$$
\Bigl \{ (c_k) \in \prod_{k \in \mC/Y} hom(F^*(X''), F^*(Y'')) \Bigl | hom(h'^*, \op{id}_{F^*(Y'')})(c_k)=hom(\op{id}_{F^*(X')}, h^*)(c_g), h:k \to g \Bigr \}
$$
which is, since $k=g \circ h$ as noted above, written as 
$$
\Bigl \{ (c_g) \in \prod_{g \in \mC/Y} hom(F^*(X'), F^*(Y')) \, \, \Bigl | \, \, c_{g \circ h} \circ h'^* =h^* \circ c_g, h:Y'' \to Y' \}.
$$
Which is nothing but $\bB^{coop}F^*(X \xrightarrow f Y)$.
\end{proof}
\begin{rem}In the above theorem we consider the category $\mS$ of sets and maps. If we want each bivariant element $c$ to be a collection of monoids or homomorphisms, then we consider the category $\mathscr Mono$ of monoids or the category $\mathscr Ab$ of abelian groups. If we want it to be a collection of $0-0$ maps of monoids or abelian groups, then we consider the category $\mathscr Mono_{0-0}$ monoids and $0-0$ maps or the catggegory $\mathscr Ab_{0-0}$ of abelian groups and $0-0$ maps. 
\end{rem}
\begin{pro}\label{coopF^*}
Let $F^*$ be a contravariant functor on a category $\mathcal C$ with the final object $pt$. We denote the set of all cohomology operations of $F^*$ by $coopF^*$. 
Then we have
$$coopF^*= \bB^{coop}F^*(pt \to pt).$$
\end{pro}
\begin{proof}
Using the notion of end, we have $coopF^* = \int_{X \in \mathcal C} hom\Bigl(F^*(X), F^*(X) \Bigr ).
$
Since $pt$ is the final object, $\mathcal C$ is the same as the slice category $\mathcal C/pt$. Let us consider the following obvious fiber square 
\begin{equation}\label{ob-square}
\xymatrix@R=20pt@C=20pt{
X \ar[d]_{\op{id}_X} \ar[r]^{p_X} & pt \ar[d]^{\op{id}_{pt}}\\
X \ar[r]_{p_X} & pt. 
} 
\end{equation}
Then we have
\begin{align*}
coopF^* & = \int_{X \in \mathcal C} hom\Bigl(F^*(X), F^*(X) \Bigr ).\\
& = \int_{p_X \in \mathcal C/pt} hom\Bigl(F^*(dom(p_X^*\op{id}_{pt})), F^*(dom(p_X)) \Bigr ) \\
& = \bB^{coop}F^*(pt \xrightarrow {\op{id}_{pt}} pt).
\end{align*}
Here $p_X^*\op{id}_{pt} = \op{id}_X$ as shown in the fiber square (\ref{ob-square}), thus $dom(p_X^*\op{id}_{pt}) = dom(\op{id}_X)=X$.
\end{proof}
\begin{defn} Using the 
pullback $p_X^*: \bB^{coop}F^*(pt \to pt) \to \bB^{coop}F^*(X  \xrightarrow {\op{id}_X} X)$ we define
\begin{align*}
coopF^*(X \xrightarrow {\op{id}_X} X) & := \op{Image} (p_X^*) \\
& = p_X^* \Bigl (\bB^{coop}F^*(pt \to pt) \Bigr ) \\
& = \Bigl \{ \left (\theta_{X'}:F^*(X') \to F^*(X') \right )_{g \in \mathcal C/X} \, \, | \, \, \theta \in coopF^* \Bigr \}.
\end{align*}
Here $X'=dom(g)$ is the domain of $g:X' \to X (\in \mathcal C/X)$.
\end{defn}
\begin{lem}\label{coopF^*in} If there is a morphism $pt \to X$, the above bivariant pullback is injective:
$$p_X^*: \bB^{coop}F^*(pt \xrightarrow {\op{id}_{pt}} pt) \lhook\joinrel\longrightarrow \bB^{coop}F^*(X \xrightarrow {\op{id}_X} X).$$
\end{lem}
\begin{proof}
The injectivity of the pullback $p_X^*$ follows from observing that, using the following fiber squares 
\begin{equation*}
\xymatrix@R=20pt@C=20pt{
pt \ar[d]_{\op{id}_{pt}} \ar[r]^{e_{pt}} & X \ar[d]_{\op{id}_X} \ar[r]^{p_X} & pt \ar[d]^{\op{id}_{pt}}\\
pt \ar[r]_{e_{pt}} & X \ar[r]_{p_X} & pt
} 
\end{equation*}
 $(\op{id}_{pt})^*= (p_X \circ e_{pt})^*= e_{pt}^* \circ p_X^*$  
 is the identity.
 \end{proof}
 
\begin{rem}
 Note that $\theta=(\theta_X)_{X \in \mathcal C}$ is indexed by the objects of $\mathcal C$, but each element in $coopF^*(X \xrightarrow {\op{id}_X} X)$ is indexed by the objects in $\mathcal C/X$. In this sense $\theta_{X'}$ can be denoted by $\theta_g$ and these two notations are interchangeably used sometimes, although we note that $\theta_g = \theta_{g'}$ if the domains of $g$ and $g'$ are the same.
\end{rem}
\begin{rem} A cohomology operation $\theta=(\theta_X)_{X \in \mC} \in coopF^*$ is one defined \emph{``globally" on the whole category $\mC$}, thus $coopF^*(X \xrightarrow {\op{id}_X} X)$ consist of such global cohomology operations \emph{``restricted"} to the slice category $\mC/X$. Thus $coopF^*(X \xrightarrow {\op{id}_X} X) \subset \bB^{coop}F^*(X \xrightarrow {\op{id}_X} X)$. On the other hand, an element $c =(c_g)_{g \in \mC/X}\in \bB^{coop}F^*(X \xrightarrow {\op{id}_X} X)$ is, from the very beginning, a cohomology operation defined \emph{``partially"} or \emph{``locally"} on the slice category $\mC/X$, therefore it requires \emph{fewer restrictions} of compatibility with the pullbacks $f^*:F^*(Y) \to F^*(X)$ for maps $f:X \to Y$, compared with a global cohomology operation $\theta =(\theta_X)_{X \in \mC}$. Thus it is quite natural that we have a strict inequality $coopF^*(X \xrightarrow {\op{id}_X} X) \subsetneqq \bB^{coop}F^*(X \xrightarrow {\op{id}_X} X)$. In other words, if the equality holds, then it means that any \emph{``locally defined"} cohomology operation $c=(c_g)_{g \in \mC/X}$ can be extended to a \emph{``globally" defined} cohomology operation $\theta=(\theta_X)_{X \in \mC}$. There are examples for which the strict inequality holds and examples for which the equality holds. In this sense, one could ask under what conditions (on $F^*$, and/or the category $\mC$ or the slice category $\mC/X$) the equality holds. However, in the rest of the paper it does not matter whether the equality holds or not. Thus, we do not discuss furthermore on this or give such examples, since these require tedious computation.
\end{rem}

\section{``Restricted" generalized cohomology operations}
In this section we discuss more about maps $c_g: F^*(X') \to F^*(Y')$ given in Definition \ref{coop-biv}.

First we consider the following definition:
\begin{defn}\label{gco} Let $F^*$ be a contravariant functor. 
An assignment $\widehat{\theta}$ defining a map $\widehat{\theta}_f: F^*(X) \to F^*(Y)$ for each continuous map  $f:X \to Y$ is called \emph{a generalized cohomology operation\footnote{This is a ``generalized" cohomology operation, not a ``generalized cohomology" operation.}} if 
for the following fiber square on the left, the following diagram on the right commutes:
\begin{equation}\label{fib.squ}
\xymatrix@R=15pt@C=20pt{
X' \ar[d]_{f'}  \ar[r]^{g'} & X \ar[d]^f \\
Y' \ar[r]_g & Y, 
} \qquad 
\xymatrix@R=15pt@C=20pt
{F^*(X') \ar[d]_{\theta_{f'}} && F^*(X)\ar[ll]_{g'^*} \ar[d]^{\theta_f}\\
F^*(Y') && F^*(Y) \ar[ll]^{g^*}.
}
\end{equation}
\end{defn}



We observed before that 
\begin{equation*}
coopF^*= \int_{X \in \mathcal C} hom\Bigl (F^*(X), F^*(X) \Bigr ).
\end{equation*}
On the other hand, the collection of all generalized cohomology operations above, denoted by $ge.coopF^*$, is described as follows, using the end:
\begin{equation}\label{ge-coop}
ge.coopF^*= \int_{f \in Mor(\mathcal C)} hom\Bigl (F^*(dom(f)), F^*(codom(f)) \Bigr ).
\end{equation}
Here $codom(f)$ is the target space $Y$ of a map $f:X \to Y$. In this sense we have
\begin{equation}\label{coop}
coopF^*= \int_{\op{id}_X \in Mor(\mathcal C)} hom \Bigl (F^*(dom(\op{id}_X)), F^*(codom(\op{id}_X)) \Bigr ).
\end{equation}
Thus the definition of a generalized cohomology operation clearly generalizes the notion of cohomology operation. However, by comparing these (\ref{ge-coop}) and (\ref{coop}) above, we can see that a generalized cohomology operation has a very \emph{strong} requirement imposed.
As far as we know, 
there seems to be no study of generalized cohomology operations. 

A generalized cohomology operation $\widehat \theta = (\widehat \theta_f)_{f \in Mor(\mathcal C}$ is indexed by the morphisms of $\mathcal C$. On the other hand, a bivariant element $c=(c_g)_{g \in \mathcal C/Y} \in \bB^{coop}F^*(X \xrightarrow f Y)$ is \emph{a ``restricted" generalized cohomology operation} in the following sense:
\begin{enumerate}
\item $f:X \to Y$ is \emph{fixed}, 
\item consider \emph{all the maps $f':X' \to Y'$ obtained by pulling back $f:X \to Y$ by any map $g:Y' \to Y$},
\item (\ref{fib.squ}) of Definition \ref{gco} is required to hold.
\end{enumerate}
 In \cite[Proposition 3.17, p.43]{Qui} (cf. \cite[\S 4 Steenrod Operations]{EG}), for the complex cobordism $U^*$ he considers 
\begin{equation}\label{pu-quillen}
P: U^{-2q}(X) \to U^{-2qk}(B \times X),
\end{equation}
which he calls the ``Steenrod $k$-th power operation". It is defined by the following composition\footnote{This kind of Steenrod power operation can be defined for \emph{any} map $f:Z \to X$, since $P:U^{-2q}(X) \to U^{-2qk}(Z)$ can be defined as the composition
$U^{-2q}(X) \xrightarrow {f^*} U^{-2q}(Z) \xrightarrow {\theta_Z} U^{-2qk}(Z)$.}
$$U^{-2q}(X) \xrightarrow {\pi^*} U^{-2q}(B \times X) \xrightarrow {\theta_{B \times X}} U^{-2qk}(B \times X)$$
where $\pi: B \times X \to X$ is the projection to the second factor and $\theta_{B \times X}$ is the Steenrod operation. 
The above assignment $P$ clearly does not satisfy (\ref{fib.squ}) in Definition \ref{gco}. 

However, \emph{if we consider a section $s:X \to B \times X$} defined by $s(x)=(\be(x),x)$ with $\be: X \to B$ is a continuous map\footnote{Namely, $\be= p_1 \circ s$, where $p_1:B \times X \to B$ is the projection to the first factor.}, then 
$$U^{-2q}(B \times X) \xrightarrow {s^*} U^{-2q}(X) \xrightarrow {\theta_X} U^{-2qk}(X)$$
becomes a ``restricted" generalized cohomology operation. 
Indeed, consider the following fiber square:
\begin{equation*}
\xymatrix@R=15pt@C=10pt
{B \times X' \ar[d]^{\pi'} \ar[rr]^{\op{id}_B \times g} && B \times X  \ar[d]_{\pi}\\
X' \ar@/^1pc/[u]^{s'} \ar[rr]_{g} && X \ar@/_1pc/[u]_{s}.
}
\end{equation*}

Then we get the following commutative diagrams:
\begin{equation*}
\xymatrix@R=15pt@C=20pt
{U^{-2q}(B \times X') \ar[d]_{(s')^*} && U^{-2q}(B \times X) \ar[ll]_{(\op{id}_B \times g)^*} \ar[d]^{s^*} \\
U^{-2q}(X') \ar[d]_{\theta_{X'}}  && U^{-2q}(X)  \ar[d]^{\theta_X}  \ar[ll]^{g^*}\\
U^{-2qk}(X')   && U^{-2qk}(X) \ar[ll]^{g^*}.
}
\end{equation*}
Hence, if we let us denote $S_{\pi}:=\theta_{X} \circ s^*$ and $S_{\pi'}:=\theta_{X'} \circ (s')^*$, 
then we get the commutative diagram
\begin{equation*}
\xymatrix@R=15pt@C=20pt
{U^{-2q}(B \times X') \ar[d]_{S_{\pi'}}  && U^{-2q}(B \times X)  \ar[d]^{S_{\pi}}  \ar[ll]_{(\op{id}_B \times g)^*} \\
U^{-2qk}(X')   && U^{-2qk}(X) \ar[ll]^{g^*}
}
\end{equation*}
Therefore the assignment $S$ does become a ``restricted" generalized cohomology operation for the projection $\pi:B \times X \to X$.

Here we note that in general $S$ does depend on the choice of a section $s:X \to B \times X$, however we also note that some different sections may give rise to the same assignment. For example, if $s_1, s_2:X \to B \times X$ are homotopic, then these two sections give rise to the same assignment. In any case, we could have many ``restricted" generalized cohomology operations for the projection $\pi:B \times X \to X$, depending on sections. 

\begin{rem}
The assinment $S$ is the pullback $s^*$ of a section followed by a cohomology operation $\theta$. But it turns out that the order of $s^*$ and $\theta$ can be exchanged. Indeed, since $\theta$ is a cohomology operation, we have the following commutative diagram
\begin{equation*}
\xymatrix@R=15pt@C=10pt
{U^{-2q}(X) \ar[d]_{\theta_X}  && U^{-2q}(B \times X)  \ar[d]^{\theta_{B \times X}}  \ar[ll]_{s^*} \\
U^{-2qk}(X)   && U^{-2qk}(B \times X) \ar[ll]^{s^*}
}
\end{equation*}
Thus we have $\theta_X \circ s^* = s^* \circ \theta_{B \times X}$. For the sake of later presentation, we will use the latter one $s^* \circ \theta_{B \times X}$, the cohomology operation followed by the pullback.
\end{rem}
\begin{defn} 
If a continuous map $f:X \to Y$ has a continuous section $s:Y \to X$, then it  shall be called a ``\emph{sectional}" map
, abusing words. 
\end{defn}

The identity map $\op{id}_X:X \to X$ and a map $a_X:X \to pt$ to a point $pt$ are trivial sectional maps with the section $s$ being $\op{id}_X$ itself and $s_x:pt \to X$  being a section choosing any point $x \in X$, respectively.

\begin{lem}\label{theta-f} Let $F^*$ be a contravariant functor and $\theta$ be a cohomology operation on $F^*$. The assignment defined by the following map $\theta_f:F^*(X) \to F^*(Y)$, for a sectional map $f:X \to Y$ with a section $s:Y \to X$, such that
$$\theta_f:= s^* \circ \theta_X: F^*(X) \xrightarrow {\theta_X} F^*(X) \xrightarrow {s^*} F^*(Y)$$
is a ``restricted" generalized cohomology operation.
\end{lem}
\begin{proof} Let us consider the following fiber squares:
\begin{equation}
\xymatrix@R=15pt@C=10pt{
Y' \ar[rr]^{g} \ar@/_1pc/[dd]_{\op{id}_{Y'}}\ar[d]^{s'}  && Y \ar[d]_s  \ar@/^1pc/[dd]^{\op{id}_Y} \\
X' \ar[d]^{f'}\ar[rr]^{g'} && X \ar[d]_f \\
Y' \ar[rr]_g  && Y.   
}
\end{equation}
Then we get the following commutative diagrams on the left, from which we get the following commutative diagram on the right:
\begin{equation*}
\xymatrix@R=15pt@C=10pt{
F^*(X') \ar[d]_{\theta_{X'}}  && F^*(X)  \ar[ll]_{g'^*} \ar[d]^{\theta_X} \\
F^*(X') \ar[d]_{s'^*}  &&F^*(X)  \ar[ll]^{g'^*}  \ar[d]^{s^*} \\
F^*(Y') && F^*(Y),  \ar[ll]^{g^*} 
}\qquad \xymatrix@R=15pt@C=10pt
{
F^*(X') \ar[d]_{\theta_{f'}}  && F^*(X)  \ar[ll]_{g'^*} \ar[d]^{\theta_f} \\
F^*(Y') && F^*(Y).  \ar[ll]^{g^*} 
}
\end{equation*}
Therefore the above assignment $\theta_f$ is a ``restricted" generalized cohomology operation.
\end{proof}
\begin{rem} The above $\theta_f$ depends on a chosen section $s:Y \to X$, thus it should be denoted, e.g., by $\theta_{f,s}$ or $\theta(s)_f$ or $\theta_f^s$, but we simply denote it as above for the sake of simplicity. If we use another section $s'$, then we denote it by $\theta'_f$.
\end{rem} 
\begin{prob}\label{problem} Given a cohomology operation $\theta$, construct a general theory of ``restricted" generalized cohomology operations $\theta_f$ for maps $f$ which are not necessarily sectional.
\end{prob}
\section{A co-operational bivaiant theory derived from cohomology operations}
In this section we construct \emph{a co-operational bivariant theory $\delta \bB^{coop}F^*(X \xrightarrow f Y)$ derived from cohomology operations}. (Here, the front symbol $\delta$ refers to ``derived".) 

From here on we assume that for a cohomology  operation $\theta$ each map  $\theta_X:F^*(X) \to F^*(X)$ is a $0-0$ map. 


\begin{defn}\label{nco-op}
$$
\delta \bB^{coop}F^*(X \xrightarrow f Y):= 
\begin{cases} 
s^* \circ coopF^*(X \xrightarrow {\op{id}_X} X), & \text{if $f$ is a sectional map with any sections $s$'s,}\\
 \bold {0} ,& \text{if $f$ is not a sectional map,}\\
\end{cases}
$$
Here $s^* \circ coopF^*(X \xrightarrow {\op{id}_X} X)$ denotes the following:
\begin{align*}
s^* \circ \, & coopF^*(X  \xrightarrow {\op{id}_X} X) := \\
&  \Bigl \{ s \circ \theta := \Bigl ((s')^* \circ \theta_{X'}:F^*(X') \xrightarrow {\theta_{X'}} F^*(X') \xrightarrow {(s')^*}  F^*(Y') \Bigr)_{g \in \mathcal C/Y} \,  \, \, \Bigl | \, \, \theta \in coopF^*, \text{sections} \,\, s \Bigr \}.
\end{align*}
Here we use the following fiber squares:
\begin{equation}\label{6-cd-yxy}
\xymatrix@R=15pt@C=20pt
{
Y' \ar[rr]^{g} \ar@/_1pc/[dd]_{\op{id}_{Y'}}\ar[d]^{s'}  && Y \ar[d]_s  \ar@/^1pc/[dd]^{\op{id}_Y} \\
X' \ar[d]^{f'}\ar[rr]^{g'} && X \ar[d]_f \\
Y' \ar[rr]_g  && Y.   
}
\end{equation}
\begin{rem} Note that the front symbol $s^*$ of $s^* \circ coopF^*(X \xrightarrow {\op{id}_X} X)$ \emph{does not} mean ``using only one section $s$", but means ``using all sections $s$'s", thus \emph{collective notation}. For $s \circ \theta := \bigl ((s')^* \circ \theta_{X'} \bigr )_{g \in \mathcal C/X}$, $s'$ for each component is the pullback of the same section $s$ by $g:X' \to X$, namely, we do not use the pullback of different sections for different components of $s \circ \theta := \bigl ((s')^* \circ \theta_{X'} \bigr )_{g \in \mathcal C/X}$.
\end{rem}
Here we emphasize that, in particular, 
$$\delta \bB^{coop}F^*(X \xrightarrow {\op{id}_X} X)= coopF^*(X \xrightarrow {\op{id}_X} X)$$ 
because the identity $\op{id}_X:X \to X$ has the unique section which is again the identity $\op{id}_X$, thus each pullbacked sections $s'$ is also the identity $\op{id}_{X'}$.


For a non-sectional map $f:X \to Y$, $\delta \bB^{coop}F^*(X \xrightarrow f Y)=\bold 0$ denotes the set of all zero maps, i.e., 
$$\bold 0 := \{0:F^*(X') \to  F^*(Y')  \, \, | \, \, g:Y' \to Y \}$$
where the zero map $0:F^*(X') \to  F^*(Y')$ sends the whole $F^*(X')$ to $0 \in F^*(Y')$. Note that the zero map is obviously a cohomology operation and also  a generalized cohomology operation. 
\end{defn}
\begin{thm}\label{coop-coho-thm} The above set $\delta \bB^{coop}F^*$ is a bivariant theory with the same co-operational bivariant operations of product, pushforwad and pullback defined in \S \ref{coop}. Hence $\delta \bB^{coop}F^*$ is a subtheory of the co-operational bivariant theory $\bB^{coop}F^*$.
\end{thm}
\begin{proof}  It suffices to show that the three co-operational bivariant operations are well-defined. 
\begin{enumerate}
\item Product: First, for the sake of convenience of reading, we recall the definition of  product: for $c  \in \delta\bB^{coop}F^*(X \xrightarrow f Y) $ and $d \in \delta\bB^{coop}F^*(Y \xrightarrow g Z) $,
$$(c \bullet d)_h:= d_h \circ c_{h'}: F^*(X') \xrightarrow {c_{h'}}  F^*(Y') \xrightarrow {d_h} F^*(Z')$$
Here we consider the following fiber squares:
\begin{equation}\label{6-cd-xyz}
\xymatrix@R=15pt@C=20pt
{X' \ar[d]_{f'} \ar[r]^{h''} & X \ar[d]^f\\
Y' \ar[d]_{g'} \ar[r]^{h'} & Y \ar[d]^g\\
Z' \ar[r]_{h} & X 
} 
\end{equation}
\begin{enumerate}
\item Let $f:X \to Y$  and $g:Y \to Z$ be both sectional maps. Then 
$$\bullet: \delta \bB^{coop}F^*(X \xrightarrow f Y) \otimes \delta\bB^{coop}F^*(Y \xrightarrow g Z) \to \delta\bB^{coop}F^*(X \xrightarrow {g \circ f} Z)$$
is well-defined. The proof of this is a bit tricky, for which we need the fact that \emph{the cohomology operation is a natural self-transformation.}
First we observe that the product $g \circ f$ is again a sectional map. Indeed, let $s_1:Y \to X$ be a section of $f$ and $s_2:Z \to Y$ be a section of $g:Y \to Z$, thus $f \circ s_1 = \op{id}_Y$ and $g \circ s_2 = \op{id}_Z$. Then we have
\begin{align*}
(g \circ f) \circ (s_1 \circ s_2) & = g \circ (f \circ s_1) \circ s_2 \\
& = g \circ \op{id}_Y \circ s_2\\
& = g \circ s_2\\
& = \op{id}_Z
\end{align*} 
Hence the product $s_1 \circ s_2$ is a section of the product $g \circ f$. Hence $g \circ f$ is also a sectional map.
Now let $c =\widetilde \theta \in \delta \bB^{coop}F^*(X \xrightarrow f Y)$ and $d= \widetilde \psi \in  \delta \bB^{coop}F^*(Y \xrightarrow g Z)$ with $\theta, \psi \in coopF^*$ and consider the following fiber squares:
\begin{equation}\label{cd-s-h-g}
\xymatrix@R=15pt@C=20pt
{
X' \ar[rr]^{h''} \ar[d]^{f'}  && X \ar[d]_f  \\
Y' \ar[d]^{g'} \ar@/^1pc/[u]^{s'_1} \ar[rr]^{h'} && Y \ar[d]_g \ar@/_1pc/[u]_{s_1}\\
Z' \ar[rr]_h \ar@/^1pc/[u]^{s'_2} && Z  \ar@/_1pc/[u]_{s_2}
} 
\end{equation}
Then by the definition we have 
$$\widetilde \theta_{h'} = s_1'^* \circ \theta_{X'}: F^*(X') \xrightarrow {\theta_{X'}} F^*(X') \xrightarrow {s_1'^*} F^*(Y'), $$
$$\widetilde \psi _h = s_2'^* \circ \psi_{Y'}: F^*(Y') \xrightarrow {\psi_{Y'}} F^*(Y') \xrightarrow {s_2'^*} F^*(Z').$$
Hence we have
\begin{equation}\label{eq-c-d}
(\widetilde \theta \bullet \widetilde \psi)_h = \widetilde \psi _h \circ \widetilde \theta_{h'} = (s_2'^* \circ \psi_{Y'}) \circ (s_1'^* \circ \theta_{X'}) = s_2'^*  \circ (\psi_{Y'} \circ s_1'^*)  \circ \theta_{X'}. 
\end{equation}
Since $\psi$ is a cohomology operation, we have the following commutative diagram:
\begin{equation}
\xymatrix@R=15pt@C=20pt
{
F^*(Y') \ar[d]_{\psi_{Y'}} && F^*(X') \ar[ll]_{s_1'^*} \ar[d] ^{\psi_{X'}}\\
F^*(Y') && F^*(X') \ar[ll]^{s_1'^*}.
}
\end{equation}
Thus we have
\begin{align*}
(\widetilde \theta \bullet \widetilde \psi)_h & = s_2'^*  \circ ( s_1'^* \circ \psi_{X'})  \circ \theta_{X'} \\
& = (s_2'^*  \circ s_1'^*) \circ (\psi_{X'}  \circ \theta_{X'})\\
& = (s_1' \circ s_2')^* \circ (\psi \circ \theta)_{X'}
\end{align*}
Since $s_1' \circ s_2'$ is a section of a sectional map $g' \circ f':X' \to Z$, which is the pullback of the sectional map $g \circ f: X \to Z$, we have
\begin{align*}
\{ (\widetilde \theta \bullet \widetilde \psi)_h  = (s_1' \circ s_2')^* \circ (\psi \circ \theta)_{X'}: & F^*(X') \to F^*(Z') \, | \, h:Z' \to Z \} \\
& \in \delta\bB^{coop}F^*(X \xrightarrow {g \circ f} Z)
\end{align*}
with $\psi \circ \theta \in coopF^*$. Therefore we have that $\widetilde \theta \bullet \widetilde \psi \in \delta\bB^{coop}F^*(X \xrightarrow {g \circ f} Z)$.
\item Suppose that either $f:X \to Y$ or $g:Y \to Z$ is not a sectional map. Then 
$$\bullet: \delta \bB^{coop}F^*(X \xrightarrow f Y) \otimes \delta \bB^{coop}F^*(Y \xrightarrow g Z) \to \delta \bB^{coop}F^*(X \xrightarrow {g \circ f} Z)$$
is well-defined. Indeed, either $\delta \bB^{coop}F^*(X \xrightarrow f Y) =\bold 0$ or $\delta \bB^{coop}F^*(Y \xrightarrow g Z) = \bold 0$, i.e., for $c \in \delta \bB^{coop}F^*(X \xrightarrow f Y)$ and $d \in \delta \bB^{coop}F^*(Y \xrightarrow g Z)$, either $d_h=0$ for any $h$ or $c_{h'} =0$ for any $h'$. If $d_h=0$ for any $h$, it is obvious that $(c \bullet d)_h = d_h \circ c_{h'} = 0$ for any $h$. If $c_{h'} =0$ for any $h'$, then it follows from ``zero-zero" condition required on a map, i.e., a map is required to send $0$ to $0$, that $(c \bullet d)_h = d_h \circ c_{h'} = 0$. Hence the product $c \bullet d =0$ belongs to 
$\delta \bB^{coop}F^*(X \xrightarrow {g \circ f} Z)$ whether $g \circ f$ is a sectional map or not, i.e., it does not matter whether the product $g \circ f$ is a sectional map or a non-sectional map.
\end{enumerate}
\item Pushforward: for $X \xrightarrow f Y \xrightarrow g Z$ 
$$f_*: \delta \bB^{coop}F^*(X \xrightarrow {g \circ f} Z) \to \delta \bB^{coop}F^*(Y \xrightarrow g Z)$$
is defined by, for $c \in \delta \bB^{coop}F^*(X \xrightarrow {g \circ f} Z)$ 
$$(f_*c)_h := c_h \circ f'^*: F^*(Y') \xrightarrow{f'^*} F^*(X')  \xrightarrow{c_h} F^*(Z').$$
Here we use the above commutative diagram (\ref{6-cd-xyz}).
\begin{enumerate}
\item If $g \circ f$ is not a sectional map, then $c_h =0$ is the zero map, thus $(f_*c)_h := c_h \circ f'^*=0$. Hence it is well-defined whether $g:Y \to Z$ is sectional or non-sectional.
\item If $g \circ f$ is a sectional map, then $g:Y \to Z$ is also a sectional map. Indeed, if there exists a section $s:Z \to X$ of $g\circ f:X \to Z$, i.e., $(g \circ f) \circ s= \op{id}_Z$, which implies that $g \circ (f\circ s)=\op{id}_Z$. Thus $f \circ s: Z \to Y$ is a section of $g:Y \to Z$. Therefore $g$ is a sectional map.
\begin{align*}
(f_*\widetilde \theta )_h & = \widetilde \theta_h \circ f'^*,  \quad \theta  \in coopF^*\\
& =(s'^* \circ \theta_{X'}) \circ f'^*\\
& = s'^* \circ (\theta_{X'} \circ f'^*) \\
& = s'^* \circ (f'^* \circ \theta_{Y'}) \quad \text{(since $\theta$ is a cohomology operation)}\\
& = (s'^* \circ f'^*) \circ \theta_{Y'}\\
& = (f' \circ s')^* \circ \theta_{Y'}:F^*(Y') \to F^*(Z')
\end{align*}
Therefore $f_*\widetilde \theta$ belongs to $\delta \bB^{coop}F^*(Y \xrightarrow g Z)$.
Here we consider the following diagrams:
\begin{equation}\label{cd-gfs}
\xymatrix@R=15pt@C=20pt{
X' \ar[rr]^{h''} \ar[d]^{f'}  && X \ar[d]_f  \\
Y' \ar[d]^{g'}\ar[rr]^{h'} && Y \ar[d]_g \\
Z' \ar[rr]_h \ar@/^1pc/[uu]^{s'} && Z  \ar@/_1pc/[uu]_{s}  
}
\quad
\xymatrix@R=15pt@C=20pt{
Y' \ar[d]^{g'}\ar[rr]^{h'} && Y \ar[d]_g \\
Z' \ar[rr]_h \ar@/^1pc/[u]^{f' \circ s'} && Z  \ar@/_1pc/[u]_{f \circ s}  
}
\quad 
\xymatrix@R=15pt@C=20pt
{
F^*(X') \ar[d]_{\theta_{X'}} & F^*(Y') \ar[l]_{f'^*} \ar[d] ^{\theta_{Y'}}\\
F^*(X') & F^*(Y') \ar[l]^{f'^*}.
}
\end{equation}
\end{enumerate}
\item Pullback:For a fiber square 
$$
\xymatrix@R=15pt@C=20pt{
X' \ar[d]_{f'} \ar[r]^{g'} & X \ar[d]^f \\
Y' \ar[r]_g & Y
}
$$
$$g^*: \delta \bB^{coop}F^*(X \xrightarrow f Y ) \to \delta \bB^{coop}F^*(X' \xrightarrow {f'} Y' )$$
is defined by, for $c \in \delta \bB^{coop}F^*(X \xrightarrow f Y) $
$$(g^*c)_h := c_{g \circ h}: F^*(X'') \to F^*(Y'').$$
Here we use the following commutative diagrams:
\begin{equation}
\xymatrix@R=15pt@C=20pt{
X'' \ar[d]_{f''} \ar[r]^{h'} & X' \ar[d]_{f'} \ar[r]^{g'} & X \ar[d]^f \\
Y'' \ar[r]_h & Y' \ar[r]_g & Y
}
\end{equation}
\begin{enumerate}
\item If $f:X \to Y$ is not a sectional map, then $\delta \bB^{coop}F^*(X \xrightarrow f Y ) =\bold 0$, hence $(g^*c)_h=0$. Thus $g^*: \delta \bB^{coop}F^*(X \xrightarrow f Y ) \to \delta \bB^{coop}F^*(X' \xrightarrow {f'} Y' )$ is well-defined whether $f':X' \to Y'$ is a sectional map or not.
\item Let $f:X \to Y$ be a sectional map with a section $s$ and $c=\widetilde \theta \in \delta \bB^{coop}F^*(X \xrightarrow f Y )$ with $\theta \in coopF^*$. Then, 
$$(g^*\widetilde \theta)_h = \widetilde \theta _{g \circ h}  = s''^* \circ \theta_{X''}: F^*(X'') \xrightarrow {\theta_{X''}} F^*(X'') \xrightarrow {s''^*} F^*(Y'').$$
Since $s''$ is also a section of $f''$ which is the pullback of $f':X' \to Y'$, we have
$$ \{ (g^*\widetilde \theta)_h = s''^* \circ \theta_{X''} \, | \, h:Y'' \to Y' \} \in \delta \bB^{coop}F^*(X' \xrightarrow {f'} Y' ).$$
Therefor the pullback is well-defined. Here we use the following commutative diagrams:
\begin{equation}
\xymatrix{
X'' \ar[rr]^{h'} \ar[d]_{f''} && X' \ar[rr]^{g'} \ar[d]_{f'} && X \ar[d]_f\\
Y'' \ar[rr]_{h} \ar@/_1pc/[u]_{s''} && Y' \ar[rr]_{g} \ar@/_1pc/[u]_{s'}  && Y \ar@/_1pc/[u]_s
}
\end{equation}
\end{enumerate}
\end{enumerate}
\end{proof}
Let $\bB$ be a bivariant theory. For a map $f:X \to Y$ we clearly have the following inclusion
$$f_*(\bB(X \xrightarrow f Y) ) \subset \bB(Y \xrightarrow {\op{id}_Y} Y)$$
by the definition of the pushforward, because $f= \op{id}_Y \circ f$, thus $\bB(X \xrightarrow f Y)= \bB(X \xrightarrow {\op{id}_Y \circ f} Y)$. Therefore, for a co-operational bivariant theory $\bB^{coop}F^*$ associated to a contravariant functor $F^*$, we do have 
$$f_*(\bB^{coop}F^*(X \xrightarrow f Y)) \subset \bB^{coop}F^*(Y \xrightarrow {\op{id}_Y} Y).$$
When it comes to the above naive co-operational bivariant theory, if $f:X \to Y$ is a sectional map, then the above inclusion becomes the equality as follows:
\begin{cor} Let $\delta\bB^{coop}F^*$ be the above naive co-operational bivariant theory. Then we have the following:
$$
f_*(\delta \bB^{coop}F^*(X \xrightarrow f Y)) = 
\begin{cases} 
\delta\bB^{coop}F^*(Y \xrightarrow {\op{id}_Y} Y)= coopF^*(Y \xrightarrow {\op{id}_Y} Y), & \text{if $f$ is a sectional map,}\\
 \bold {0}  \, ( \in \delta \bB^{coop}F^*(Y \xrightarrow {\op{id}_Y} Y)\, ) ,& \text{if $f$ is not a sectional map,}\\
\end{cases}
$$
\end{cor}
\begin{proof}
The case when $f:X \to Y$ is not a sectional map is obvious. Let $f:X \to Y$ be a sectional map with a section $s$. Consider the above commutative diagram (\ref{6-cd-yxy}). Let us take an element 
$$\{s'^* \circ \theta_{X'}:F^*(X') \xrightarrow {\theta_{X'}} F^*(X') \xrightarrow {s'^*}  F^*(Y')  \, \, | \, \, g:Y' \to Y \} \in s^* \circ coopF^*(X \xrightarrow {\op{id}_X} X).$$
Then it follows from the definition of the co-operational bivariant pushforward that
\begin{align*}
f_*(s'^* \circ \theta_{X'}) & = (s'^* \circ \theta_{X'}) \circ f'^*  \quad \text{(by the definition of the pushforward)} \\
& = s'^* \circ (\theta_{X'} \circ f'^*) \\
& = s'^* \circ (f'^* \circ \theta_{Y'}) \quad \text{(since $\theta$ is a cohomology operation)}\\
& = (s'^* \circ f'^*) \circ \theta_{Y'}\\
& = (f' \circ s')^* \circ \theta_{Y'}\\
& = (\op{id}_{Y'})^* \circ \theta_{Y'}\\
& = \theta_{Y'}.
\end{align*}
Therefore we obtain $f_*(\delta \bB^{coop}F^*(X \xrightarrow f Y)) = \delta \bB^{coop}F^*(Y \xrightarrow {\op{id}_Y} Y)$.
\end{proof}
Let $T:F^* \to G^*$ be a natural transformation between two contravariant functors. Suppose that for a certain set $coop_TF^*$ of cohomology operations of $F^*$ the natural transformation $T:F^* \to G^*$ is \emph{compatible with} the cohomology operations $coop_TF^*$ of $F^*$ and the cohomology operations $coopG^*$ of $G^*$, which is denoted by 
$T: coop_TF^* \to coopG^*,$
which shall be also called ``a natural transformation of cohomology operations", abusing words. 
This means the following (cf. Definition \ref{4-defn-bcoop} and (\ref{cd-7})). Let $\theta$ be a cohomology operation of $F^*$ and let $T(\theta)=:\theta^T$. Then for any map $f:X \to Y$ the following cubic diagram commutes:
\begin{equation}
\xymatrix@R=15pt@C=10pt
{ F^*(Y) \ar[dd]_{\theta_Y } \ar[rd]^{f^*} \ar[rr]^{T} && G^*(Y) \ar'[d][dd]^(.4){(\theta^T)_Y} \ar[rd]^{f^*} \\
& F^*(X) \ar[dd]_(.3){ \theta_X} \ar[rr]^(.3){T}  &&  G^*(X) \ar[dd]^{(\theta^T)_X} \\
F^*(Y) \ar'[r] [rr]_{T \quad \quad }  \ar[rd]_{f^*} && G^*(Y) \ar[rd]_{f^*} \\
& F^*(X) \ar[rr] _{T }  &&  G^*(X).
}
\end{equation}
\begin{ex}
We consider the Chern character $ch: \op{Vect}(-) \to H^*(-)\otimes \mathbb Q$ 
for an example of $T:F^* \to G^*$. For that purpose, we consider the following ``polynomial operation".

For a vector bundle $E$ we define
$$E^{\oplus k}:= \underbrace{E \oplus E \oplus \cdots \oplus E}_k, \quad E^{\otimes k}:= \underbrace{E \otimes E \otimes \cdots \otimes E}_k.$$
For $n=0$, $E^0$ is defined to be the trivial line bundle $\jeden$.
Let $p(x) \in \mathbb Z_{\geq 0}[x]$ be an integral polynomial with non-negative coefficients, say $p(x) = \sum_{i=0}^n a_ix^i$, then we define\footnote{Here we mention the notion of \emph{finite vector bundle} due to Andr\'e Weil \cite{Weil} (cf. \emph{essentially finite vector bundle} \cite{Nori1, Sz, BV}), although we do not discuss it in this paper. A vector bundle $E$ is called \emph{finite} if there exist two different polynomials $p(x)$ and $q(x)$ with non-negative coefficients such that $p(E) \cong q(E).$
A trivial vector bundle is finite, since, for any trivial vector bundle $\jeden^r$ of rank $r$, consider a constant polynomial $p(x)=r$ and $q(x)=x$ ($p(x)\not = q(x)$), then we do have $p(\jeden^r)=\jeden^r = q(\jeden^r).$ 
}
$$p(E) := \bigoplus_{i=0}^n(E^{\otimes i})^{\oplus a_i}.$$ 
which shall be denoted by $p(x)(E)$. 
For a polynomial $p(x) \in \mathbb Z_{\geq 0}[x]$, 
we have the 
commutative diagram for any continuous map $f:X \to Y$:
\begin{equation*}
\xymatrix@R=20pt@C=10pt{
\op{Vect}(X) \ar[d]_{p(x)}  && \op{Vect}(Y)  \ar[ll]_{f^*} \ar[d]^{p(x)} \\
\op{Vect}(X) && \op{Vect}(Y).  \ar[ll]^{f^*} 
}
\end{equation*}
i.e., $p(x)(f^*E) = f^*\left (p(x)(E) \right )$, which is easy to see.
For the cohomology functor $H^*(-)$ and an element $\alp \in H^*(X)$ we denote
$$\alp^n := \underbrace{\alp \cup \alp \cup \cdots \cup \alp}_n.$$
Then, similarly, 
we define 
$p(x)(\alp):=p(\alp) = \sum_{i=0}^n a_i\alp^i$
 and for $f:X \to Y$ we have
$$p(x)(f^*\alp) = p(f^*\alp) = \sum_{i=0}^n a_i(f^*\alp)^i = \sum_{i=0}^n f^*(a_i\alp^i) = f^*(\sum_{i=0}^n a_i\alp^i) = f^*\left (p(x)(\alp) \right ).$$
The Chern character $ch: \op{Vect}(-) \to H^*(-)$ commutes with these polynomial operations, i.e., 
$$p(x)(ch(E)) = ch (p(x)(E)),$$
which follows from $ch(E \oplus F) =ch(E) +ch(F)$ and $ch(E\otimes F) =ch(E) \cup ch(F).$
Therefore both of the following diagrams commute:
\begin{equation}\label{ch-diagram}
\xymatrix@R=20pt@C=10pt
{
H^*(X) \ar[d]_{p(x)}  && H^*(Y)  \ar[ll]_{f^*} \ar[d]^{p(x)} \\
H^*(X) && H^*(Y).  \ar[ll]^{f^*} 
}\qquad 
\xymatrix@R=15pt@C=10pt
{
\op{Vect}(X) \ar[d]_{p(x)}  \ar[rr]^{ch} && H^*(Y) \ar[d]^{p(x)} \\
\op{Vect}(X) \ar[rr]_{ch} && H_*(X). 
}
\end{equation}
Now, let $coop_{ch}\op{Vect}$ be the set of all the polynomial cohomology operations defined above, thus for $p(x) \in coop_{ch}\op{Vect}$ we have $p(x)(E):=p(E)$. Then the above $T(\theta)=:\theta^T$ means that $p(x)^{ch}:=p(x)$, i.e., a polynomial cohomology operation of $H^*(-)\otimes \mathbb Q$ defined before, i.e., $p(x)(\alp):=p(\alp)$. Thus the Chern character $ch:\op{Vect}(-) \to H^*(-) \otimes \mathbb Q$ gives rise to a natural transformation
$$ch: coop_{ch}\op{Vect}(-) \to coopH^*(-) \otimes \mathbb Q,$$
namely we have the following cubic commutative diagram:
\begin{equation*}
\xymatrix@R=15pt@C=10pt
{ 
\op{Vect}(Y) \ar[dd]_{p(x) } \ar[rd]^{f^*} \ar[rr]^{ch} && H^*(Y)\otimes \mathbb Q \ar'[d][dd]^(.4){p(x)} \ar[rd]^{f^*} \\
& \op{Vect}(X) \ar[dd]_(.3){ p(x)} \ar[rr]^(.3){ch}  &&  H^*(X)\otimes \mathbb Q  \ar[dd]^{p(x)} \\
\op{Vect}(Y) \ar'[r] [rr]_{ch \quad \quad }  \ar[rd]_{f^*} && H^*(Y)\otimes \mathbb Q \ar[rd]_{f^*} \\
& \op{Vect}(X) \ar[rr] _{ch }  &&  H^*(X)\otimes \mathbb Q .
}
\end{equation*}
\end{ex}
\begin{ex}
Let $c\ell: \op{Vect}(-) \to H^*(-)\otimes \mathbb Q$ be any multiplicative characteristic class such as Chern class, Todd class, L-class, i.e., $c\ell(E \oplus F) = c\ell(E) \cup c\ell(F)$.
Then the commutative diagram on the right in (\ref{ch-diagram}) is replaced by the following:
\begin{equation*}\label{diagram-ch}
\xymatrix@R=20pt@C=10pt
{
\op{Vect}(X) \ar[d]_{\ell_a(x)}  \ar[rr]^{c\ell} && H^*(Y) \ar[d]^{\mu_a(x)} \\
\op{Vect}(X) \ar[rr]_{c\ell} && H_*(X),
}
\end{equation*}
because for any $\ell_a(x)=ax \in \mathbb Z_{\geq 0}[x]$ and $\mu_a(x)=x^a \in \mathbb Z_{\geq 0}[x]$ we have
$$c\ell(\ell_a(x)(E)) = c\ell(\ell_a(E))=c\ell(E^{\oplus a})= c\ell(E)^a = \mu_a(c\ell(E)) = \mu_a(x)(c\ell(E)).$$
Now, let $coop_{c\ell}\op{Vect}$ be the set of all cohomology operations defined by a linear one $\ell_a(x)=ax$. 
Then the above $T(\theta)=:\theta^T$ means that $c\ell (\ell_a(x))= \ell_a(x)^{c\ell}:=\mu_a(x)$ where $\mu_a(x)=x^a$. 
Thus a natural transformation $c\ell:\op{Vect}(-) \to H^*(-) \otimes \mathbb Q$ gives rise to a natural transformation
$$c\ell: coop_{c\ell}\op{Vect}(-) \to coopH^*(-) \otimes \mathbb Q,$$
namely we get the following cubic commutative diagram:
\begin{equation*}
\xymatrix@R=15pt@C=10pt
{ 
\op{Vect}(Y) \ar[dd]_{\ell_a(x) } \ar[rd]^{f^*} \ar[rr]^{c\ell} && H^*(Y)\otimes \mathbb Q \ar'[d][dd]^(.4){\mu_a(x)} \ar[rd]^{f^*} \\
& \op{Vect}(X) \ar[dd]_(.3){ \ell_a(x)} \ar[rr]^(.3){c\ell}  &&  H^*(X)\otimes \mathbb Q  \ar[dd]^{\mu_a(x)} \\
\op{Vect}(Y) \ar'[r] [rr]_{c\ell \quad \quad }  \ar[rd]_{f^*} && H^*(Y)\otimes \mathbb Q \ar[rd]_{f^*} \\
& \op{Vect}(X) \ar[rr] _{c\ell }  &&  H^*(X)\otimes \mathbb Q .
}
\end{equation*}
\end{ex}

\begin{cor} Let $T:F^* \to G^*$ be a natural transformation between two contravariant functors and we have a natural transformation $T: coop_TF^* \to coopG^*$ as above. Then $T: coop_TF^* \to coopG^*$ can be extended to a Grothendieck transformation between the naive co-operational bivariant theories:
$$\delta T:\delta\bB_T^{coop}F^*(X \to Y)  \to \delta\bB^{coop}G^*(X \to Y),$$
provided that the three properties listed in Corollary \ref{cor-nat} hold.
Here $\delta\bB_T^{coop}F^*(X \to Y)$ is defined by replacing $coopF^*(Y \xrightarrow {\op{id}_Y} Y)$ by  $coop_TF^*(Y \xrightarrow {\op{id}_Y} Y)$ in Definition \ref{nco-op}.
\end{cor}
\begin{proof} Let the situation be as above. Then we have the following commutative diagrams:
 \begin{equation}\label{cd-FG}
\xymatrix@R=20pt@C=10pt{
F^*(X') \ar[r]^{T} \ar[d]_{\theta_{X'}}  & G^*(X')\ar[d]^{\theta^T_{X'}}  \\
F^*(X')  \ar[d]_{s'^*}\ar[r]^{T} & G^*(X') \ar[d]^{s'^*} \\
F^*(Y')  \ar[r]_T & G^*(Y')  
}
\end{equation}
Then we define $\delta T:\delta\bB_T^{coop}F^*(X \to Y)  \to \delta\bB^{coop}G^*(X \to Y)$ as follows:
\begin{enumerate}
\item If $f: X \to Y$ is not a sectional map, $\delta T(\bold 0):=\bold 0.$
\item  If $f: X \to Y$ is a sectional map, then  
$$\delta T(\{s'^* \circ \theta_{X'}\, \,  | \, \, g:Y' \to Y \} ):=\{s'^*  \circ \theta^T_{X'}\, \,  | \, \, g:Y' \to Y  \},$$
for which see the above diagrams (\ref{cd-FG}).
\end{enumerate}
\end{proof}
\noindent
{\bf Acknowledgements}:
Some work of this paper was done while the author was staying at The Erwin Schr\"odinger International Institute for Mathematics and Physics (ESI) and attended the two-week Workshop ``Algebraicity and Transcendence for Singular Differential Equations" (Oct. 7 -- Oct.19, 2024). The author would like to express sincere thanks to one of the organizers, Professor Herwig Hauser (Universit\"at Wien), for the invitation to the workshop and the staff of ESI for providing the author with a very nice atmosphere to work in. The author also would like to express his sincere gratitude to the referee for a report with constructive suggestions and comments.  
The research is supported by JSPS KAKENHI Grant Number JP23K03117.

\end{document}